\documentclass[12pt]{amsart}
\usepackage[left=1in,top=1in,right=1in,bottom=1in]{geometry}
\usepackage{mathrsfs}
\usepackage[english]{babel}
\usepackage{amsmath,amsthm,amsfonts,amssymb,epsfig,xcolor, hyperref}
\usepackage{bbm,booktabs,float,mathtools,siunitx,tikz}
\usepackage[shortlabels]{enumitem}
\usepackage{breqn}
\usepackage{stmaryrd}
\usepackage{relsize}

\def\P{\ensuremath\mathbb{P}}

\def\E{\ensuremath\mathbb{E}}

\numberwithin{equation}{section}
\newtheorem{thm}{Theorem}
\newtheorem{lem}[thm]{Lemma}
\newtheorem{prop}[thm]{Proposition}
\newtheorem{cor}[thm]{Corollary}

\newtheorem*{rem}{Remark}

\numberwithin{thm}{section}

\newcommand{\PP}{\mathbf{P}}

\newcommand{\re}{{\rm Re}}

\begin{document}
	\title{Upper bounds on large deviations of Dirichlet $L$-functions in the $q$-aspect}
	
	\author{Louis-Pierre Arguin}
	
	\address{L.P. Arguin, Mathematical Institute, University of Oxford, UK}
	\email{louis-pierre.arguin@maths.ox.ac.uk}
	
	\author{Nathan Creighton}
	\address{N. Creighton, Mathematical Institute, University of Oxford, UK}
	\email{creighton@maths.ox.ac.uk}
	
	\begin{abstract}
		We prove a result on the large deviations of the central values of even primitive Dirichlet $L$-functions with a given modulus.	For $V\sim \alpha\log\log q$ with $0<\alpha<1$, we show that \begin{equation}\nonumber\frac{1}{\varphi(q)} \# \left\{\chi \text{ even, primitive mod }q: \log \left|L\left(\chi,\frac{1}{2}\right)\right| >V\right\}\ll  \frac{e^{-\frac{V^2}{\log\log q}}}{\sqrt{\log\log q}}.\end{equation}
		This yields the sharp upper bound for the fractional moments of central  values of Dirichlet $L$-functions proved by Gao, upon noting that the number of even, primitive characters with modulus $q$ is $\frac{\varphi(q)}{2}+O(1).$
		The proof is an adaptation to the $q$-aspect of the recursive scheme developed by Arguin, Bourgade and Radziwi\l\l{ }for the local maxima of the Riemann zeta function, and applied by Arguin and Bailey to the large deviations in the $t$-aspect. We go further and get bounds on the case where $V=o(\log\log q)$. These bounds are not expected to be sharp, but the discrepancy from the Central Limit Theorem estimate grows very slowly with $q$. The method involves a formula for the twisted mollified second moment of central values of Dirichlet $L$-functions, building on the work of Iwaniec and Sarnak. 	\end{abstract}
	\date{May 31, 2024}
	\maketitle

	\section{Introduction}
	This paper is on the central value of Dirichlet $L$-functions with  fixed modulus $q$. The most famous result in the study of the distribution of the values of Dirichlet $L$-functions is Selberg's Central Limit Theorem. In \cite{Sel}, it is shown that the normalised values of the logarithm of the Riemann zeta function in a dyadic interval high up the critical line have a Gaussian limiting distribution. More precisely, for any measurable set $E\subset \mathbb C,$
	\begin{equation}\label{Selberg}
		\lim_{T\to \infty} \frac{1}{T}{\rm meas}\left\{ t\in [T,2T]: \frac{\log\left(\zeta\left(\frac{1}{2}+it\right)\right)}{\sqrt{\frac{1}{2}\log\log T}}\in E\right\}=\frac{1}{2\pi}\int_Ee^{-\frac{x^2+y^2}{2}}dxdy.
	\end{equation}
A remarkable feature of Selberg's Central Limit Theorem is that it gives information about the distribution of the imaginary part of the logarithm of the Riemann zeta function on the critical line, the fluctuation of which is connected to zeroes of the Riemann zeta function, without assuming the Riemann Hypothesis.
  
	Much work has been done to generalise Equation \eqref{Selberg}, both in extending the range of the Gaussian behaviour, and proving central limit theorems for other problems related to Dirichlet $L$-functions. For example, in \cite{PW}, it is shown that the real parts of the logarithm of distinct Dirichlet $L$-functions in a dyadic interval high up the critical line behave like independent Gaussian variables. More precisely, for any integer $N$, if $(\chi_j)^N_{j=1}$ is a sequence of distinct Dirichlet characters, then   Theorem 1.3 of \cite{PW} states  that for $t$  uniformly distributed in the interval $[T,2T]$ with $T$ large,  the random vector $$\left(\log\left|L\left(\frac{1}{2}+it,\chi_1\right)\right|,\dots \log\left| L\left(\frac{1}{2}+it,\chi_N\right)\right|\right)$$ approximately has the distribution of a real Gaussian random variable with mean $\mathbf{0}$ and variance $\frac{1}{2}\left(\log\log T\right) \mathbf{I}_N.$
	
	The most famous conjecture about the non-vanishing of Dirichlet $L$-functions is the Generalised Riemann Hypothesis, which predicts that all non-trivial zeroes of primitive Dirichlet $L$-functions lie on the critical line with real part $\frac{1}{2}.$ However, even on the critical line, one can consider whether it is possible for any Dirichlet $L$-function to vanish at a certain point. It is conjectured that no primitive Dirichlet $L$-function vanishes at the central point. In \cite{Pratt}, it is shown that for large values of $q$, at least $50.073\% $ of the central values of $L$-functions associated to primitive characters of modulus $q$ are non-vanishing. In this article we are able to show the large deviations of order $\log\log q$ are bounded above by a Gaussian tail. If one were able to extend this to show that the large deviations beyond the order $\sqrt{\log q\log\log q}$ were also bounded by the same Gaussian tail, one would be able to show the non-vanishing of all central values, however the authors were unable to extend their results to this range.  
	
	There are many similarities between the distribution in the $q$-aspect of central values of Dirichlet $L$-functions with a large modulus, and the distribution in the $t$-aspect of the values of a given Dirichlet $L$-function, such as the Riemann zeta function, high up the critical line. The analysis of the distribution of values of $L$-functions via their moments is a widely investigated area, and parallels can be drawn here between the $q$-aspect and the $t$-aspect.  The Keating-Snaith Conjecture from \cite{KS}
	predicts asymptotics for the moments of the Riemann zeta function. It  conjectures that
	
	\begin{equation}
		\lim_{T\to\infty} \frac{1}{T(\log T)^{\lambda^2}}\int^T_0 \left|\zeta\left(\frac{1}{2}+it\right)\right|^{2\lambda}dt= C_\lambda,\end{equation}
	where $C_\lambda$ is an explicit constant.
	For the $q$-aspect, an analogous result is conjectured in \cite{RS}, for all $\lambda>0:$ \begin{equation}\label{qconj}
		\lim_{q\to\infty} \frac{1}{\varphi(q)(\log q)^{\lambda^2}}\sum_{\chi \text{ primitive mod }q} \left|L\left(\chi,\frac{1}{2}\right)\right|^{2\lambda}= C_{1,\lambda},\end{equation}
	where $C_{1,\lambda}$ is an explicit constant.
	
	 In \cite{gao}, it is shown that
	\begin{equation}\label{gaomoment}
		\frac{1}{\varphi(q)}\mathlarger{\sum}_{\chi \text{ primitive mod $q$}} \left|L\left(\chi,\frac{1}{2}\right)\right|^{2\beta}\asymp (\log q)^{\beta^2},
	\end{equation} for all $\beta\in (0,1)$. This means that, up to a constant, these moments of $\left|L\left(\chi,\frac{1}{2}\right)\right|$, for $\chi$ ranging over the primitive characters with modulus $q$, match the moment generating function of a Gaussian random variable with variance $\frac{1}{2}\log\log q$.
	Moreover, we have a lower bound for all the integer moments of the conjectured correct order from \cite{RS}.
	Here, it is shown that, for $k$ any fixed natural number and $q$ a large prime,
	\begin{equation}\label{lowermoment}
		\frac{1}{\varphi(q)}\sum_{\chi \text{ primitive mod $q$}}\left|L\left(\chi,\frac{1}{2}\right)\right|^{2k}\gg_k (\log q)^{k^2}. 
	\end{equation}
	
	Equation \eqref{lowermoment} is proven in \cite{RS} for the case where $q$ is a large prime, but the method of proof easily generalises to the general case for large values of $q$ not congruent to $2$ modulo $4.$	
	This means the integer moments of the central values are bounded below by the moment generating function of a Gaussian random variable. We prove further that the large deviations of the real part of the logarithm of the central value are bounded above by a Gaussian tail. 
	
	In \cite{AB}, a large deviations result is proven for the Riemann zeta function; this paper shows an analogous result holds for the central values of Dirichlet $L$-functions with a given modulus. It was necessary to prove certain results used in the proof of the large deviations result in \cite{AB} for the context of central values of Dirichlet $L$-functions; most notably, the twisted mollifier formula for second moments of the Riemann zeta function needed to be adapted to the $q$-aspect, cf.~Theorem \ref{twistthm}. As in \cite{Sarnak}, we restrict our attention to look at just the even primitive characters; the case of odd primitive characters yields the same results. 
	\subsection{Main results}
	\begin{thm} \label{LargeDev} Let $q$ be a large natural number, with $q\ne 2 \text{ mod } 4.$ Suppose that \\$V\sim \alpha\log\log q,$ with $0<\alpha<1$.
		Then 
		\begin{equation}\frac{1}{\varphi(q)} \# \left\{\chi \text{ even,  primitive mod }q: \log \left|L\left(\chi,\frac{1}{2}\right)\right| >V\right\}\ll  \frac{e^{-\frac{V^2}{\log\log q}}}{\sqrt{\log\log q}}.\end{equation}
	\end{thm}	
	In Section \ref{CorSec}, we deduce an upper bound on the fractional moments of the central values of primitive Dirichlet $L$-functions, in line with Gao's result, Equation \eqref{gaomoment}.
	\begin{cor}\label{momentcor}
		Let $0<\beta<1$ and $q$ be a large natural number, with $q\ne 2$ mod $4$. Then we have \begin{equation}
			\frac{1}{\varphi(q)}\mathlarger{\sum}_{\chi\text{ even,  primitive mod }q } \left|L\left(\chi,\frac{1}{2}\right)\right|^{2\beta}\ll  (\log q)^{\beta^2},
		\end{equation}
		where the  implicit constant is uniform for $\beta$ lying in any interval $(0,B)$, for any $B<1$.
	\end{cor}
	Note that the range for $\beta$ in Corollary \ref{momentcor} is controlled by the range of $\alpha$ for which Theorem \ref{LargeDev} is valid. Indeed, if we were able to extend the range of Theorem \ref{LargeDev} to be valid for $0<\alpha < A$ for some $A>1$, then we would be able to extend our bound on the moments in Corollary \ref{momentcor} to the range $0<\beta<A.$
	
	In \cite{AB}, they re-prove a result of Heap, Radziwi\l\l{ }and Soundararajan \cite{HRS19} which is the analogous result to Corollary \ref{momentcor} in the $t$-aspect.
	
	\begin{cor}[Corollary 1.2 in \cite{AB}]
		Let $0<\beta<2$. Then there exists a constant $C_\beta$ such that
		\begin{equation}\label{blowup}
			\frac{1}{T}\int^{2T}_T \left|\zeta(1/2+it)\right|^{2\beta} \le C_\beta (\log T)^{\beta^2},
		\end{equation}
		for all $T$ sufficiently large.
	\end{cor}	
	The bounds in \cite{AB} have $C_\beta\to \infty$ as $\beta\to 0$. However, by tightening the bounds when we adapt their method, we are able to avoid our constant in the $q$-aspect becoming unbounded near $\beta=0$.
	Because we are able to tighten our bounds, we are able to attain (weaker) bounds which include the case when $V=o(\log\log q).$
	\begin{thm}\label{largerange}
		Let $q$ be a large natural number, with $q\ne 2$ mod $4$ and $V>0$ with $V=o(\log\log q)$. Let $\mathscr{L}$ be as in Equation \eqref{L}. Then \begin{equation}\frac{1}{\varphi(q)}\left\{\chi \text{ even,  primitive mod }q: \log \left|L\left(\chi,\frac{1}{2}\right)\right| >V\right\}\ll  \frac{\mathscr{L}e^{-\frac{V^2}{\log\log q}}}{\sqrt{\log\log q}}.\end{equation}
	\end{thm}
	Note that $\mathscr{L}$ grows very slowly in $q$. Indeed, for a given natural number $n$, the maximum value of $q$ for which  $\mathscr{L}=n$ grows like the Ackermann function $F(4,n).$
	
	\subsection{Notation and structure of the paper}
	We will be taking expectations and sums over different sets of characters with modulus $q$.
	\begin{itemize}
		\item The set of primitive even characters with modulus $q$, denoted $+$.
		\item The set of all even characters with modulus $q$, denoted $\oplus$.
		\item The set of all primitive characters with modulus $q$, denoted $\ast$.
	\end{itemize}
	Where we take probabilities, $\P$ will always denote the probability where a character $\mathfrak{X}$ is drawn uniformly at random from the space of primitive even characters with modulus $q$. 
	
	
	\subsection{Structure of the proof}  
	Theorems \ref{LargeDev} and \ref{largerange} are proved in Section \ref{devsec}. The proof of Corollary \ref{momentcor} follows easily and the details are given in Section \ref{CorSec}.
	Theorem \ref{LargeDev} is derived from the recursive scheme first developed for the $t$-aspect in \cite{ABR20}. The idea is to approximate $\log | L(1/2, \chi)|$ by Dirichlet polynomials of increasing length and whose values are restricted. 
	To this end,  we split up the primes contributing most to the Dirichlet $L$-functions into intervals. We iteratively define steps by setting $q_0=1.5,$ and setting \begin{equation} q_l=\exp\left(\frac{\log q}{\log_{l+1}(q)^{\mathbf{s}}}\right)\end{equation} for $l\ge 1,$ where $\log_l$ means the natural logarithm iterated $l$ times. Here, $\mathbf s$ is a large constant, and taking $\mathbf s>10^5$ as in \eqref{eqn: s parameter}. 
	We then set \begin{equation}\label{defn}
		n_l=\log\log q_l,\quad l\geq 1,
	\end{equation}
	which turns out to be the suitable scale for considering the variance of Dirichlet polynomials supported on primes up to $q_l.$
	We halt our steps at $\mathscr{L}$, the largest value of $l$ such that 
	\begin{equation}\label{L}
		\exp(10^6 (\log\log q-n_l)^{10^5}e^{n_{l+1}})=q^{\frac{10^6}{(\log_{l+2}q)^{10^5-\mathbf s}}}
		\le q^{\frac{1}{100}}.
	\end{equation}
	
	For $k>0$, we define a truncated sum of the formal logarithm of the Dirichlet $L$-series at $\frac{1}{2}$:
	\begin{equation}\label{tildes}
		\tilde{S}_k=\mathlarger{\sum}_{p\le e^{e^k}}\frac{\chi(p)}{p^{\frac{1}{2}}}+\frac{\chi(p)^2}{2p}.
	\end{equation}
	The absolute value of the $L$-function at $\frac{1}{2}$ may be controlled by  the real part of the above series.
	We further define
	\begin{equation}
		\label{eqn: S}
		S_k=\mathlarger{\sum}_{p\le e^{e^k}}\Re\left(\frac{\chi(p)}{p^{\frac{1}{2}}}+ \frac{\chi(p^2)}{2p}\right).
	\end{equation}	
	We expect for $k$ large, $S_k$ to be close to $\Re(\log L\left(\chi,\frac{1}{2}\right))$.
	Indeed, we only consider the values of $S_k$ for $k$ close to $\log\log q$, at the time steps $k=n_l$, for $1\le l\le \mathscr{L}$.
	
	We also need to control the difference $\log \left|L(\chi,\frac{1}{2})\right|- S_k$.
	To do so, we define a mollifier of the contribution of the primes in the interval $(q_{l-1},q_l]$.
	Given $1\le l\le \mathscr{L}$, set
	\begin{equation} \label{lmol}
		M_l(\chi,s)=\mathlarger{\sum}_{\substack{p|n\implies p\in (q_{l-1},q_l]\\ \Omega_l(n)\le 10(n_l-n_{l-1})^{10^5}}} \frac{\mu(n)\chi(n)}{n^{s}},
	\end{equation}
	where $\Omega_l(n)$ denotes the number of prime factors with multiplicity of $n$ in the interval $(q_{l-1},q_l]$. The idea is that the mollifiers are long enough that $M_1...M_l\left(\chi,\frac{1}{2}\right)$ should mollify all primes up to $q_l$. We want to show $\left|M_1...M_l\right|$ to be a good approximation of $\exp(-S_{n_l})$. Then, we will have $\left|M_1...M_lL(\chi,\frac{1}{2})\right|$ not too small or large, so we will have mollified the $L$-function at $\frac{1}{2}$ successfully.
	
	The mollifier described above is part of a larger family of Dirichlet polynomials which vary corresponding to the even primitive characters mod $q$ by a scaling of their coefficients.	Namely, if $\chi$ is an even primitive character mod $q$, then its associated Dirichlet $L$-function is \begin{equation}L(\chi,s)=\mathlarger{\sum}_n \frac{\chi(n)}{n^{s}}.\end{equation}
	For a sequence of real coefficients $(a_m)$, which is uniform for all even primitive  Dirichlet characters with modulus $q$ and supported on the integers $1\le m\le M$, we consider the associated Dirichlet polynomial \begin{equation} \label{qcoef}Q(\chi,s)=\mathlarger{\sum}_n \frac{a_n\chi(n)}{n^{s}}.\end{equation}
	We say $Q=Q(\chi)$ is a {\it degree $e^{n_l}$ well-factorable polynomial} if we can write \begin{equation}\label{Qsplits}
		Q=\prod^l_{\lambda=1} Q_\lambda,\end{equation} 
	where $Q_\lambda$ is a Dirichlet polynomial with the support of the coefficients restricted by their prime factors, and with the coefficients having the same dependence on $\chi$ as in \eqref{qcoef},
	say 
	\begin{equation}
		Q_\lambda(\chi,s)= \mathlarger{\sum}_{\substack{p|m\implies p\in (q_{\lambda-1},q_{\lambda}]\\ \Omega_\lambda(m) \le 10(n_\lambda-n_{\lambda-1})^{10^4}}}\frac{\gamma_m \chi(m)}{m^{s}},\end{equation}
	i.e., each factor $Q_\lambda(\chi,s)$ has coefficients supported on integers with all prime factors in the interval $(q_{\lambda-1},q_\lambda]$, and  is not too long.
	We are only interested in the central values of the Dirichlet polynomials; for a Dirichlet polynomial $Q$ with such a property we will write $Q$ for $Q\left(\chi,\frac{1}{2}\right).$ 
	
	We show that if $Q$ is a degree $e^{n_l}$ well-factorable Dirichlet polynomial, then $M_1...M_l$ mollifies $L$ successfully, after a twist by $Q$. We consider the Euler product of $L$, which doesn't converge at $\frac{1}{2}$, but can still give good intuition.
	Since $M_1...M_l$ is a good mollifier of $L$, we would expect only the primes greater than $q_l$ to affect the value of $M_1...M_lL$. Hence, if $Q$ is a Dirichlet polynomial which only depends on the primes up to $q_l$, then we would expect $LM_1...M_l$ and $Q$ to be weakly dependent. 
	This is a crucial step of the recursive scheme.
	We show this in Section \ref{Twistsec}, where we prove the following theorem:
	
	\begin{thm}\label{twistthm}
		\textit{Let $Q$ be a degree $e^{n_l}$ well-factorable polynomial and $L$ and the mollifiers $M_1,...,M_l$ be as defined above. Then the mollified twisted second moment satisfies } 
		\begin{equation}\E_{+}[\left|LM_1...M_lQ\right|^2]\ll  \frac{\log q}{\log q_l}\E_{\oplus}[\left|Q\right|^2].
		\end{equation}
	\end{thm}
	
	The expectation over $+$ is the expectation where $\chi$ is drawn uniformly at random from the primitive even characters mod $q$, whilst the expectation over $\oplus$ is where $\chi$ is drawn uniformly from all even characters mod $q$ (so includes the imprimitive even characters), and has simpler orthogonality relations. These relations are discussed in Section \ref{Momentsec}.
	We see that the ratio of the expectations in Theorem \ref{twistthm} is bounded by a factor independent of the twist $Q$, which suggests $LM_1...M_l$ and $Q$ are weakly dependent.
	We need a slightly stronger mollified twisted moment than in Theorem $1.5$, to handle twists by polynomial functions of the real parts of Dirichlet polynomials, such as $S_k-S_{k-1}.$
	
	We require a similar version to being well-factorable in the real case.
	We say a function $F=F(\chi)$ is {\it $l$-sufficient} if we may write:
	\begin{equation}\label{lsuff}
		F=\prod^l_{j=1} F_j,
	\end{equation}
	where for each $1\le j\le l$, we have
	\begin{equation}F_j=K_j{\Big(}\Re{\Big(}\sum_{m} \chi(m) b^{j}_{(m)}{ \Big)}{\Big)},	\end{equation}
	for some sequence of complex numbers $b^{j}_{(m)}$ and polynomial $K_j$ of degree $d_j$ with real coefficients.
	The coefficient  $b^{j}_{(m)}$ are assumed to be $0$ unless $p|m\implies p\in (q_{j-1},q_j]$ and $\Omega_j(m) \le \nu_j$ for $\nu_j= \frac{10(n_j-n_{j-1})^{10^4}}{d_j}$.
	The bound on the number of prime factors ensures that the length of $F_j$ is less than $10(n_j-n_{j-1})^{10^4}$.
	In Section \ref{Twistsec} we will prove:
	\begin{thm}\label{realtwistthm} \textit{Let $1\le l\le \mathscr{L}$ and  let $F$ be $l$-sufficient. Then the  twisted second moment mollified by $F$ satisfies } \begin{equation}\E_{+}[\left|LM_1...M_lF\right|^2]\ll  \frac{\log q}{\log q_l}\E_{\oplus}[\left|F\right|^2].\end{equation}
	\end{thm}
	By applying Theorem \ref{realtwistthm} for suitable choices of $F,$ we will deduce Theorem \ref{twistthm}.
	In order to prove Theorem \ref{realtwistthm}, we require bounds on the twisted moments of Dirichlet $L$-functions by real parts of Dirichlet polynomials.
	We will follow the formulae for twists by the whole Dirichlet polynomial (not just the real part), i.e.,
	$\E_{+}\left[\left|LM\right|^2\right]$,
	where $M$ is a short Dirichlet polynomial. In \cite{Bui}, the authors give a bound where $M$ has length up to $q^{\frac{51}{101}}$. However, for our purposes it suffices to follow the proof of the simpler formula in \cite{Sarnak}, for Dirichlet polynomials of length smaller than $q^{\frac{1}{2}}.$\\
	
	\noindent {\bf Acknowledgements.}
	Both authors thank Jon Keating for insightful discussions on the subject. The research of L.-P. A. is partially supported by the NSF grant DMS 2153803. The research of N.C. is supported by the EPRSC grant EP/W524311/1.

	\section[A]{Moments in the {$q$}-aspect}\label{Momentsec}
	In this section we collate some simple bounds on the moments of even Dirichlet polynomials with a fixed modulus, which will be needed in the proof of Theorems \ref{LargeDev}, \ref{twistthm} and \ref{realtwistthm}.
	\subsection{Moments over all even characters with a fixed modulus}
	The following identity is exact, as opposed to the $t$-aspect, where the diagonal terms only give a leading order approximation.
	\begin{lem}
		Let $q$ be a modulus, $N< \frac{q}{2}$ a positive integer, and $(a_n)^N_{n=1}$ be a sequence of complex numbers. Then \begin{equation}\label{orthog}
			\E_{\oplus} \left[\left|\mathlarger{\sum}^N_{n=1} a_n\chi(n)\right|^2\right]=\mathlarger{\sum}^N_{n=1}\left |a(n)\right|^2,
		\end{equation} where the expectation is taken over even characters mod $q$.
	\end{lem}
	\begin{proof} Given integers $1\le m,n\le N,$ we have \begin{equation}\mathlarger{\sum}\limits_{\chi \text{ mod $q$ even}}\chi(n)\overline{\chi(m)}=\\
			\mathlarger{\sum}\limits_{\chi \text{ mod $q$}}\left(\frac{\chi(n)+\chi(-n)}{2}\right)\overline{\chi(m)}.\end{equation}
		This equals $\frac{\phi(q)}{2}$ if $n=m$ and $0$ otherwise, since we cannot have $n\equiv -m \text{ mod $q$}$.
		Upon substitution into the right-hand side of \eqref{orthog}, we get: 
		\begin{equation}
			\E_{\oplus}\left[ \left|\mathlarger{\sum}^N_{n=1} a_n\chi(n)\right|^2\right]=\frac{2}{\varphi(q)}\mathlarger{\sum}\limits_{1\le n,m\le N} a_n\overline{a_m} \mathlarger{\sum}\limits_{\chi \text{ mod $q$ even}} \chi(n)\overline{\chi(m)}.\end{equation}
		Only the diagonal terms survive, and we conclude. \end{proof}
	We deduce the following splitting theorem:
	\begin{lem}\label{Splitform} If $A,B$ are two Dirichlet polynomials varying with $\chi$, each of length shorter than $  \sqrt{\frac{q}{2}}$ with support on integers with prime factors in disjoint sets, then \begin{equation}\E_{\oplus}\left[\left|AB(1/2,\chi)\right|^2\right]=\E_{\oplus}[|A(1/2,\chi)|^2]E_{\oplus}\left[\left|B(1/2,\chi)\right|^2\right].\end{equation}\
	\end{lem} 
	Note again that the expectation splits exactly due to the exact orthogonality relations, which is a stronger result than in the $t$-aspect (cf \cite{AB}, Lemma 13).
	We also require results on the moments of real parts of Dirichlet polynomials in the $q$-aspect. It is convenient to introduce the random variables $(X(p), p \text{ primes})$ that are independent and uniformly distributed on the unit circle. 
	The random variable $X(n)$ is then defined multiplicatively by
	\begin{equation}
		\label{eqn: X}
		X(n)= X(p_1)^{a_1} \dots X(p_r)^{a_r},
	\end{equation} 
	for an integer $n$ with prime factorisation $n=p_1^{a_1}...p_r^{a_r}$.
	\begin{lem}\label{Realexp}
		Let $q$ be large, $K$ be a polynomial of degree d with real coefficients, $N< \frac{q^{\frac{1}{2d}}}{2}$ be a positive integer, $B=(b_n)^N_{n=1}$ be a sequence of real variables, with $b_n=0$ unless $n$ is a prime or a square of a prime.
		Then $\E_{\oplus}\left[K\left(\Re \left(\sum^N_{n=1} b_n\chi(n)\right)\right)^2\right]=\E \left[K\left(\Re\left(\sum^N_{n=1} b_n X(n)\right)\right)^2\right]$. 
	\end{lem}
	\begin{proof}
		
		Since $K$ is a polynomial with real coefficients, we may write \begin{equation}\label{Krep}
			K\left(\Re \left(\sum^N_{n=1} b_n\chi(n)\right)\right)=\mathlarger{\sum}\limits_{\substack{j,k\le N^{d}\\ p|jk\implies p\le N}} c_{j,k,B,K}\chi(j)\overline{\chi(k)}\end{equation} for certain coefficients $c_{j,k,B,K}$ depending on the polynomial $K$ and the coefficients $B=(b_n)^N_{n=1},$ satisfying \begin{equation}\label{symmetry}
			c_{j,k,B,K}=\overline{c_{k,j,B,K}}.
		\end{equation}
		We may cancel any common factors of $k$ and $j$ using the relation: $\chi(n)\overline{\chi(n)}=1$ for any character with modulus $q$ and any integer $n$ coprime to $q$.
		For coprime integers $m$ and $n$ with at most $N^d$ with prime factors all at most $N$, let \begin{equation}\label{Cdef}
			C_{m,r,B,K}=\sqrt{mr}\mathlarger{\sum}\limits_{\substack{\frac{j}{(j,k)}=m\\
					\frac{k}{(j,k)}=r}} c_{j,k,B,K}.
		\end{equation}
		Note that the sum is zero unless $m$ and $r$ are coprime, with all their prime factors at most $N$ so we may drop these conditions. By \eqref{symmetry}, we see
		\begin{equation}\label{Csymmetry}
			C_{m,r,B,K}=\overline{C_{r,m,B,K}}.
		\end{equation}	
		Here, the factor $\sqrt{mr}$ in \eqref{Cdef} was introduced to give the same scaling convention for coefficients at the point $s=\frac{1}{2}$ used in other sections of the article, such as for the Dirichlet polynomial \eqref{qcoef}.
		From \eqref{Krep} we see,		
		\begin{equation}\label{Kcop}K\left(\Re \left(\sum^N_{n=1} b_n\chi(n)\right)\right)=\mathlarger{\sum}\limits_{m,r\le N^d\\ 
			} \frac{ C_{m,r,B,K}}{\sqrt{mr}} \chi(m)\overline{\chi(r)}.\end{equation}
		Hence \begin{equation}\label{Ksq}K\left(\Re \left(\sum^N_{n=1} b_n\chi(n)\right)\right)^2=\mathlarger{\sum}\limits_{m_1,m_2, r_1,r_2\le N^d\\
			}\chi(m_1m_2)\overline{\chi(r_1r_2)} \frac{C_{m_1,r_1,B,K} C_{m_2,r_2,B,K}}{\sqrt{m_1r_1m_2r_2}}.
		\end{equation}
		The restriction $N<\frac{q^{\frac{1}{2d}}}{2}$ ensures that we cannot have $m_1m_2\equiv -r_1r_2\text{ mod $q$, }$and that\\
		\begin{equation}\label{realdiag}
			m_1m_2\equiv r_1r_2\text{ mod $q$ }\iff m_1m_2=r_1r_2.
		\end{equation}
		Moreover, the restrictions on the support of the coefficients $(m_1,r_1)=(m_2,r_2)=1$ ensure that \eqref{realdiag} is only achieved in the support if \begin{equation}
			m_1=r_2, m_2=r_1.
		\end{equation}
		Using the relation \eqref{orthog} for the Dirichlet polynomial in \eqref{Ksq}, we see:
		\begin{equation}\label{thetalink}
			\E_{\oplus}\left[K\left(\Re \left(\sum^N_{n=1} b_n\chi(n)\right)\right)^2\right]=\mathlarger{\sum}\limits_{m,r\le N^d} \frac{C_{m,r,B,K}C_{r,m,B,K}}{mr}.
		\end{equation}
		Finally, applying \eqref{Csymmetry}, we see 
		\begin{equation}
			\E_{\oplus}\left[K\left(\Re \left(\sum^N_{n=1} b_n\chi(n)\right)\right)^2\right]=\mathlarger{\sum}\limits_{m,r\le N^d} \frac{\left|C_{m,r,B,K}\right|^2}{mr}.
		\end{equation}
		The corresponding identity with $\chi(n)$ replaced by $X(n)$ is proved similarly
		
		\begin{equation}\E \left[K\left(\Re\left(\sum^N_{n=1} b_n X(n)\right)\right)^2\right]=\mathlarger{\sum}\limits_{m,r\le N^d} \frac{\left|C_{m,r,B,K}\right|^2}{mr}.\end{equation}
		since the orthogonality relation \eqref{orthog} holds by the definition of $X(n)$.

		%
		%
	\end{proof}
	
	Using Lemma \ref{Realexp}, we deduce a splitting relation for the real part of Dirichlet polynomials.
	\begin{lem}\label{Realsplit}
		Let $1\le l\le \mathscr{L}$ and $F$ be $l$-sufficient as defined in \eqref{lsuff}.
		Then
		\begin{equation}
			\E_{\oplus}\left[F^2\right]=\prod^l_{j=1}\E_{\oplus}\left[K_j\left(\Re\left(\sum_{n} b^{(j)}_n \chi(n)\right)\right)^2\right].
		\end{equation}
	\end{lem}
	\begin{proof}
		Using the notation from the proof of Lemma \ref{Realexp}, we see that for each $1\le j\le l$, the coefficients $C_{m,r,B^{(j)},K_j}$ are supported on coprime integers with prime factors all lying in $(q_{j-1},q_j]$, and we similarly drop this condition from the sum. Then
		\begin{equation}
			K_j\left(\Re\left(\sum_{n} b^{(j)}_n \chi(n)\right)\right)^2=\sum_{ m_1,m_2,r_1,r_2}
			\chi(m_1m_2)\overline{\chi(r_1r_2)} \frac{C_{m_1,r_1,B^{(j)},K_j} C_{m_2,r_2,B^{(j)},K_j}}{\sqrt{m_1m_2r_1r_2}}.
		\end{equation}
		Moreover, for all values of $r$ and $m,$ we have \begin{equation}\label{Splitsym}
			C_{m,r,B^{(j)},K_j}=\overline{C_{r,m,B^{(j)},K_j}},
		\end{equation}
		and since $b_n(j)$ is zero unless $\Omega_j(n) \le \nu_j$, we see $C_{r,m,B^{(j)},K_j}$ is zero unless $\Omega_j(r),\Omega_j(m)\le 10 (n_j-n_{j-1})^{10^4}$. By Equation \eqref{orthog}, we obtain: \begin{equation}\label{realpartexp}
			\E_{\oplus}\left[\Re\left(K_j\left(\sum_{n} b^{(j)}_n \chi(n)\right)\right)^2\right]=\sum_{m,r}
			\frac{\left|C_{m,r,B^{(j)},K_j}\right|^2}{mr}.
		\end{equation}
		Given coprime integers $u$ and $v$, each having all their prime factors at most $q_l$, we may write \begin{equation}u=\prod^l_{j=1} m_j, v=\prod^l_{j=1} r_j, \end{equation} 
		where $m_j$ and $r_j$ are coprime integers with all their prime factors lying in the interval $(q_{j-1},q_j].$
		Then set \begin{equation}\label{deftc}
			\tilde{C}_{u,v}= \prod^l_{j=1} C_{m_j, r_j, B^{(j)}, K_j},  
		\end{equation}
		where we suppress the dependence of $\tilde{C}_{u,v}$ on the sequences $B^{(j)}$ and the polynomials $K_j$ for ease of notation. 
		The relation \eqref{Splitsym} means that \begin{equation}\label{Jointsym}
			\tilde{C}_{u,v}=\overline{\tilde{C}_{v,u}}.
		\end{equation}
		We note that $\tilde{C}_{u,v}$ is supported on coprime integers $u$ and $v$ with all their prime factors at most $q_l$, so we drop this from the summation criteria to show that for any even primitive character modulo $q$,
		\begin{equation}\label{jointK}\prod^l_{j=1}K_j \left(\Re\left(\sum_{n} b^{(j)}_n \chi(n)\right)\right)=\mathlarger{\sum}\limits_{u,v} \frac{\tilde{C}_{u,v}}{\sqrt{uv}} \chi(u)\overline{\chi(v)}.\end{equation}
		The restriction on the support of the coefficients $\tilde{C}_{u,v}$ means $\tilde{C}_{u,v}=0$ unless $u,v\le \frac{q^{{\frac{1}{2}}}}{2},$ so we may apply \eqref{orthog} to yield
		\begin{equation}\label{jointexpold}
			\E_{\oplus}\left[\prod^l_{j=1}K_j \left(\Re\left(\sum_nb^{(j)}_n \chi(n)\right)\right)^2\right]=\mathlarger{\sum}\limits_{u,v} \frac{\left|\tilde{C}_{u,v}\right|^2}{uv}.
		\end{equation}
		
		Splitting into the components from prime factors in each interval, and using the support of the coefficients $\tilde{C}_{u,v},$ we may write the right-hand side of \eqref{jointexpold} as
		\begin{equation}
			\mathlarger{\sum}\limits_{\forall 1\le j\le l, p|m_jn_j\implies p\in (q_{j-1},q_j]} \frac{\left|\tilde{C}_{\prod^l_{j=1}m_j, \prod^l_{j=1}n_j}\right|^2}{\prod^l_{j=1}m_jn_j}.
		\end{equation}
		Using \eqref{deftc}, we may write this as
		\begin{equation}
			\prod^l_{j=1} 	\mathlarger{\sum}\limits_{ p|m_jn_j\implies p\in (q_{j-1},q_j]} \frac{\left|{C}_{m_j,n_j, B^{(j)},K_j}\right|^2}{m_jn_j}.
		\end{equation}
		
		Finally, using \eqref{realpartexp}, we may write the above expression as 
		\begin{equation}
			\prod^l_{j=1}	\E_{\oplus}\left[K_j\left(\Re\left(\sum_{n} b^{(j)}_n \chi(n)\right)\right)^2\right] .\end{equation}
		This completes the proof of Lemma \ref{Realsplit}.
	\end{proof}
	
	\subsection{Moments over even primitive characters with a given modulus}
	In most of the proofs, we take the expectation of Dirichlet polynomials over all even characters. However, we require a bound on the expectation of certain Dirichlet polynomials where we restrict to just the even primitive characters.
	
	\begin{lem}\label{primlem}
		Let $(a_n)^N_{n=1}$
		be a sequence of  length $N\le q^{\frac{5}{6}}$.
		Then we have the following bound on the expectation over even primitive characters with modulus $q$:
		\begin{equation}\label{primorthog}
			\E_{+}\left[\left|\sum^N_{n=1}a_n \chi(n)\right|^2\right] \ll \sum^N_{n=1}\left|a_n\right|^2   
		\end{equation}
	\end{lem}
	Note that \eqref{primorthog} is an upper bound for the case where there is no restriction to primitive characters in \eqref{orthog}, for the restricted length $N\ll q^{\frac{5}{6}}.$
	\begin{proof}
		Expanding the square out in \eqref{primorthog} shows that the left-hand side is equal to
		\begin{equation}
			\frac{1}{\varphi(q)}\sum^N_{j,k=1} a_j\overline{a_k} \left(\sum^+_{\chi \text{ even,  primitive mod $q$ }} \chi(j)\overline{\chi(k)}\right)
		\end{equation}
		Using the Cauchy-Schwarz inequality, this may be bounded as
		\begin{equation}
			\ll \frac{1}{\varphi(q)}\mathlarger{\sum}\limits_{j\le N}\left|a_j\right|^2\left|\mathlarger{\sum}\limits_{k\le N}\sum^+_{\chi \text{ even,  primitive mod $q$ }}  \chi(j)\overline{\chi(k)}\right|.
		\end{equation}
		So it suffices to show that
		\begin{equation}\label{kbound}
			\frac{1}{\varphi(q)}\mathlarger{\sum}\limits_{k\le N} \left|\sum^+_{\chi \text{ even,  primitive mod $q$ }} \chi(j)\overline{\chi(k)}\right|\ll 1,
		\end{equation} uniformly for each $j$ coprime to $q$.
		
		Equation 3.2 in \cite{Sarnak} gives the sum over primitive characters with modulus $q$ (even and odd), at an individual value. Using Mobius inversion and the orthogonality of characters, they show, for any $m$ coprime to $q$:
		\begin{equation}\label{avg}
			\sum^\ast_{\chi\text{ primitive, even or odd mod $q$}} \chi(m)=\mathlarger{\sum}\limits_{\substack{vw=q\\ m\equiv 1\text{ mod }q}}\mu(v)\varphi(w). 
		\end{equation}
		We have
		\begin{equation}\label{reflect}
			\sum^*_{\chi \text{ even, primitive mod $q$}} \chi(m)=\frac{1}{2}\left(\sum^\dagger_{\chi\text{ even or odd mod $q$}} \chi(m)+\chi(-m)\right).
		\end{equation}
		
		When substituting the bounds from \eqref{reflect} and \eqref{avg}, we see that \eqref{kbound} is bounded as:
		\begin{equation}\label{Mobinvert}
			\ll \frac{1}{\varphi(q)}\left(\mathlarger{\sum}\limits_{k\le N} \left|\mathlarger{\sum}\limits_{\substack{vw=q\\ j\equiv k\text{ mod $w$}}} \mu(v)\varphi(w)\right| +\mathlarger{\sum}\limits_{k\le N} \left|\mathlarger{\sum}\limits_{\substack{vw=q\\ j\equiv -k\text{ mod }w}} \mu(v)\varphi(w)\right|\right).\end{equation}
		If we write $t=(j-k,q)$ and $u=(j+k,q)$,
		then we transform \ref{Mobinvert} into:
		
		\begin{align} \label{groupfactor}
			\ll \frac{1}{\varphi(q)}\mathlarger{\sum}\limits_{t,u|q}\left(\#\{k\le N: (j-k,q)=t\}\left|\mathlarger{\sum}\limits_{ w|t } \mu(\frac{q}{w}) \varphi(w)\right|\right.\\\nonumber+\left.\#\{k\le N: (j+k,q)=u\}\mathlarger{\sum}\limits_{w|u }\left| \mu(\frac{q}{w}) \varphi(w)\right|
			\right).
		\end{align}	
		
		But given any factors $t$ and $u$ of $q$, we have 
		\begin{equation}
			\#\{k\le N: (j-k,q)=t\}\ll \frac{N}{t}+1,\quad \#\{k\le N: (j+k,q)=u\}\ll \frac{N}{u}+1,
		\end{equation}
		so that \eqref{groupfactor} may be bounded as 
		\begin{equation}\label{doublefactor}
			\ll \frac{1}{\varphi(q)} \mathlarger{\sum}\limits_{u|q} \left(\frac{N}{u}+1\right)\left|\mathlarger{\sum}\limits_{{ w|u }} \mu(\frac{q}{w}) \varphi(w)\right|.
		\end{equation} 
		We may write the prime factorisations of $q,u$ and $w$ as:
		\begin{equation}\label{qprimes}
			q=p_1^{a_1}...p_n^{a_n},\quad u=p_1^{b_1}...p_n^{b_n},\quad w=p_1^{c_1}...p_n^{c_n},
		\end{equation}
		respectively, where the $p_i$ are distinct, $a_i\ge 1$ for all $i$, and $0\le c_i \le b_i\le a_i$. Examining the term $\mu(\frac{q}{w}),$ we see that the contribution to \eqref{doublefactor} is zero, unless $b_i\in \{a_i-1, a_i\}$ for all $i$.
		Then the choices for $c_i$ giving non-zero contribution to \eqref{doublefactor} are:
		\begin{equation}\begin{cases}
				c_i=a_i-1& b_i=a_i-1\\
				c_i\in \{a_i-1, a_i\}& b_i=a_i\end{cases}.
		\end{equation}
		If we define $S=\{1\le i\le n: b_i=a_{i}-1\}$, to represent the prime factorisation of $u$, then the expression in \eqref{doublefactor} becomes:
		
		\begin{equation}
			\frac{1}{\varphi(q)}\mathlarger{\sum}\limits_{S\subset \{1,...,n\}} 
			\left(\frac{N}{q}\prod_{i\in S} p_i+1\right) \left|\prod_{i\in S} \mu(p_i)\varphi(p_i^{a_i-1})\prod_{i\notin S} \left(\mu(1) \varphi(p_i^{a_i})+\mu(p_i)\varphi(p_i^{a_i-1})\right)\right|.  
		\end{equation}
		We split the sum up by the factors in the first bracket.
		We have to bound:
		
		\begin{equation}
			\mathscr{B}_1=	\frac{1}{\varphi(q)}	
			\mathlarger{\sum}\limits_{S\subset \{1,...,n\}}  \left|\prod_{i\in S} \mu(p_i)\varphi(p_i^{a_i-1})\prod_{i\notin S} \left(\mu(1) \varphi(p_i^{a_i})+\mu(p_i)\varphi(p_i^{a_i-1})\right)\right|,
		\end{equation}
		
		and \begin{equation}
			\mathscr{B}_2=\frac{N}{q\varphi(q)}\mathlarger{\sum}\limits_{S\subset \{1,...,n\}} 
			\prod_{i\in S} p_i \left|\prod_{i\in S} \mu(p_i)\varphi(p_i^{a_i-1})\prod_{i\notin S} \left(\mu(1) \varphi(p_i^{a_i})+\mu(p_i)\varphi(p_i^{a_i-1})\right)\right|.
		\end{equation}  
		We will show that both $	\mathscr{B}_1$ and $	\mathscr{B}_2$ are of order $1$, completing the proof of Lemma \ref{primlem}.
		
		\textbf{Bounding $\mathscr{B}_1$}
		The expression for $\mathscr{B}_1$ splits into a product over the prime factors of $q$.
		If the factorisation is as in \eqref{qprimes}, then we see:
		\begin{equation}
			\mathscr{B}_1=	
			\prod^r_{i=1}\frac{1}{\varphi(p_i^{a_i})}\left(\left| \varphi(p_i^{a_i-1})\right|+\left|\varphi(p_i^{a_i})-\varphi(p_i^{a_i-1})\right|\right).
		\end{equation}
		Since $\varphi(p_i^{a_i-1})\le \varphi(p_i^{a_i})$, we see that the inner term is identically $1$, and so $\mathscr{B}_1=1$.

		\textbf{Bounding $\mathscr{B}_2$}
		The term for $\mathscr{B}_2$ also splits as a product over prime factors of $q.$
		If the factorisation is as in \eqref{qprimes}, then we see:
		\begin{equation}\label{Nbound}
			\mathscr{B}_2\le \frac{N}{q^{\frac{5}{6}}}\prod^{n}_{i=1} \frac{1}{p_i^{\frac{a_i}{6}}\varphi(p_i^{a_i})} \left(\left| p_i \varphi(p_i)^{a_i-1}\right|+\left| \varphi(p_i^{a_i})-\varphi(p_i^{a_i-1})\right|\right).
		\end{equation}
		For $p_i> 2$, we bound the term in  the product using that
		\begin{equation}\frac{1}{p_i} \varphi(p_i^{a_i})\le \varphi(p_i^{a_i-1})\le \frac{\varphi(p_i^{a_i})}{p_i-2}.\end{equation}
		For $p_i=2$ we simply use that  \begin{equation}\left|\varphi(p_i^{a_i-1})\right| \ll \left|\varphi(p_i^{a_i})\right|,\end{equation} to show that the term associated to $p_i=2$ is bounded. 
		Under these bounds we see that  the term in \eqref{Nbound} yields:
		\begin{equation}
			\mathscr{B}_2\ll \frac{N}{q^{\frac{5}{6}}}\prod_{p_i>2} \frac{1}{p_i^{\frac{a_i}{6}}}\left(\frac{p_i}{p_i-2}+\frac{p_i-1}{p_i}\right).
		\end{equation}
		The term in the product: \begin{equation}
			\frac{1}{p_i^{\frac{a_i}{6}}}\left(\frac{p_i}{p_i-2}+\frac{p_i-1}{p_i}\right)
		\end{equation}
		can only be greater than $1$ for finitely many choices of $p_i$ and $a_i$, and so we see that 
		$\mathscr{B}_2\ll \frac{N}{q^{5/6}}$
		By the bounds on $N$, this yields $\mathscr{B}_2\ll 1$.	
	\end{proof}
	\begin{rem}
		{\rm
			From the proof, one can show that \begin{equation}\label{errorexp}\E_{+}\left[\left|\sum^N_{n=1}a_n\chi(n)^2\right|\right]\ll \left(1+O\left(\frac{N}{q^{\frac{5}{6}}}\right)\right)\sum^N_{n=1}\left|a_n\right|^2,\end{equation}
			and hence deduce a version of Lemma \ref{Splitform} for splitting expectations over primitive characters. One could then perform all the expectations in this paper over the even primitive characters with modulus $q$, without adding in all the non-primitive characters. However, we take the view that, for simplicity, it is easier to just work with all the even characters, to avoid the cumulative error term in \eqref{errorexp}.}
	\end{rem}
	Following this proof, we can also get the expectation of polynomials of real parts of Dirichlet characters over primitive even characters, as opposed to over all even characters in Lemma \ref{Realsplit}. We observe that the definition of an $l$-sufficient function means that the \textit{length} is $O(q^{\frac{1}{100}}).$
	\begin{lem}\label{realprim}
		Let $1\le l\le \mathscr{L}$ and $F$ be $l$-sufficient as defined in \eqref{lsuff}.
		Then the expectation of F over primitive even characters satisfies \begin{equation}\E_{+}\left[F^2\right]\ll\prod^l_{j=1}\E_{\oplus}\left[K_j\left(\Re\left(\sum_{n} b^{(j)}_n \chi(n)\right)\right)^2\right].\end{equation}
	\end{lem}

	\begin{lem} 
		\label{lem: partmoment}
		Given k an integer, let $\tilde{S}_k$ be as defined in Equation \eqref{tildes}. If $ n< m$, and $2n\le \frac{\log q}{e^m}$, then \begin{equation} \label{partmoment}
			\E_{+}\left[ \left|\tilde{S}_n-\tilde{S}_m\right|\right]^{2k}\ll  k!(m-n+1)^k.
		\end{equation}
	\end{lem}
	Note this is an  analogue of Lemma 3 in \cite{Sound}, but we give more details.
	We first explain some intuition behind the result. We split the sums $\tilde{S}_m$ and $\tilde{S}_n$ into the contributions from the primes and the contribution from the square of the primes. The expectation of  the contribution of the primes give the $m-n$ term.	For fixed $k$, the contribution of the square of the primes is negligible, and in \cite{Sound} is neglected.
	But to take $k$ increasing, we must consider the contribution of the squares of primes in \eqref{partmoment}. This gives an upper bound of $(m-n+1)^k,$ rather than $(m-n)^k$, which is what we would expect from just the primes.
	We first consider the contribution of the primes, without the squares. We require a lemma on the expected value of Dirichlet polynomials supported on primes.
	\begin{lem}\label{Soubound}  If $a(p)$ is a sequence of complex numbers, and $x$ is a natural number such that $x^k\le q,$ then
		\begin{equation}\mathlarger{\sum}_{\chi}^*\left(\mathlarger{\sum}\limits_{e^{e^n}\le p\le x} \left|\frac{a(p)\chi(p)}{p^{\frac{1}{2}}}\right|^{2k}\right) \ll qk!\left(\mathlarger{\sum}\limits_{e^{e^n}\le p\le x}\frac{|a(p)^2|}{p}\right)^k.\end{equation}
	\end{lem}
	\begin{proof} This follows from the proof of Lemma 3 in \cite{Sound}, using the orthogonality relations for Dirichlet characters from Lemma \ref{primlem}.\end{proof}
	
	\begin{proof}[Proof of Lemma \ref{lem: partmoment}]
		We use the same notation for expanding the $k^{th}$ moment of a Dirichlet series as in Lemma 3 in \cite{Sound}, so that if 
		$N=\prod^R_{i=1} p_i^{\alpha_i}$,
		then we set \begin{equation}
			a_{k,x}(N)={k\choose \alpha_1\ ...\ \alpha_R}\prod^R_{i=1}a(p_i)^{\alpha_i}. 	
		\end{equation}
		The contribution of the primes to the sum $\tilde{S}_n-\tilde{S}_m$ is: \begin{equation}
			U_{1}:=\mathlarger{\sum}\limits_{e^{e^n}\le p \le e^{e^m}}\frac{\chi(p)}{p^{\frac{1}{2}}}.
		\end{equation}
		Hence applying Lemma \ref{Soubound} we see \begin{equation}
			\E_{+} \left[\left|U_1\right|^{2k}\right]\ll  k!\left(\mathlarger{\sum}\limits_{e^{e^n}\le p\le e^{e^m}}\frac{1}{p}\right)^k\ll  k!(m-n)^k,\end{equation}
		where the last inequality follows by Mertens's estimate.

		We need to evaluate the contribution of the  terms from the squares of primes.
		The contribution of the squares of primes to the sum $\tilde{S}_m-\tilde{S}_n$ is: \begin{equation}U_2:=\mathlarger{\sum}\limits_{p\le x} \frac{\chi(p^2)}{2p}.
		\end{equation}
		But we see on expansion that
		\begin{equation}U_2^k =\mathlarger{\sum}\limits_{r}\frac{a_{k,x}(r)\chi(r^2)}{r}.\end{equation}
		and applying the same arguments as in the proof of Lemma 3 in \cite{Sound}, with the orthogonality relations for Dirichlet characters shows that  \begin{equation}\E_{+}\left[\left|U_2\right|^{2k}\right]\ll k! \left(\mathlarger{\sum}\limits_{{e^{e^m}}\le p\le e^{e^n}} \frac{1}{2p^2}\right)^k
			\le k! (1)^k.
		\end{equation}
		Hence  \begin{equation}\E_{+}\left[\left|  \tilde{S}_m-\tilde{S}_n\right|^{2k}\right]
			=\E_{+}\left[\left| U_1+U_2\right|^{2k}\right]\ll  k!(m-n+1)^k\end{equation} by H{\"o}lder's Inequality. 
		
	\end{proof}

	By setting $k=\lceil \frac{V^2}{m-n+1}\rceil $ in \eqref{partmoment}, we get by Stirling's formula and Markov's inequality:
	\begin{lem} Let q be a large integer, with $q\ne 2$  mod $4,$ $V>0$ and $0<n<m\le q_{\mathscr{L}}$.  Then
		\begin{equation}\label{softbound}
			\frac{1}{\varphi(q)}\#\{\chi \text{ even primitive mod }q: \left|\tilde{S}_m-\tilde{S}_n\right|>V\ll  \frac{V+1}{(m-n+1)^{\frac{1}{2}}}\exp\left(-\frac{V^2}{m}\right).
	\end{equation}\end{lem} 
	We require a stronger result for the moments of the real parts of the truncated sums.
	\begin{lem}\label{Strongmoment}
		Let $n<m<\mathscr{L}$ and $2k\le \frac{\log q}{3 \log q_m}$. Then
		\begin{equation}\label{hardS}
			\E_{+}\left[\left|S_{q_m}-S_{q_n}\right|^{2k}\right]\ll   \frac{(2k)!}{2^{2k} k!} (\log \log q_m-\log \log q_n)^{k}.\end{equation}\end{lem}
	Note that \eqref{hardS} is stronger than the analogous Lemma A.2 in \cite{AB}.   \begin{proof}
		Let
		\begin{equation}\label{achoice}
			a_{r,m,n}=\begin{cases}
				1 & p\text{ a prime in }(q_n, q_m]\\
				0 & \text{otherwise}
			\end{cases}	.
		\end{equation}
		Then set $U_{m,n}=\sum\limits_r \frac{a_{r,m,n}}{r^{1/2}} \chi(r)$, and $T_1=\Re U_{m,n}$
		to be the real part of the sum $U_{m,n}.$ Then $T_1$ is the contribution of the primes to the difference $S_m-S_n$.
		Then if we set $T_2=(S_m-S_n)-T_1$, then $T_2$ is the contribution of the square of primes to the difference.
		Clearly, we have
		\begin{equation} \E_{+}\left[T_1^{2k}\right]=\frac{1}{2^k}\sum^{k}_{j=0}{2k\choose j}\E_{+}[U_{m,n}^j \overline{U_{m,n}}^{2k-j}].\end{equation}
		Using the same orthogonality relations as in the proof of  \eqref{partmoment}, and observing that the choice of $k$ means $(S_m-S_n)^k$ has length 
		$q_m^{2k} \le (q_m)^{\frac{ \log q}{3 \log q_m}}=q^{\frac{1}{3}}$ we obtain
		\begin{equation}\label{realexp2}
			\E_{+}\left[T_1^{2k}\right]	\ll \frac{1}{2^k}\sum^k_{j=0}{2k\choose j}\E[\hat{U}_{m,n}^j \overline{\hat{U}_{m,n}}^{2k-j}],\end{equation}
		where   \begin{equation}\hat{U}_{m,n}=\sum\limits_{q_n\le p\le q_m} \frac{a_p}{p^{\frac{1}{2}}} X(p),\end{equation}
		and $X(r)$ is defined in \eqref{eqn: X}.
		Under this definition, \eqref{realexp2} becomes:
		\begin{equation}\E(T_1)^{2k}		\ll \E[\Re(\hat{U}_{m,n})^{2k}].\end{equation}	 	
		We need the following proposition, which is a straightforward computation.
		\begin{prop} \textit{ For any complex numbers $a_p$, if $X(p)$ are as defined in \eqref{eqn: X}, then we have: }\begin{equation}\label{bessapp}\E\left[\sum\limits_{x\le p\le y} (\Re  \ a_pX(p))^{2k}\right]\le\frac{(2k!)}{2^k k!}(s^2)^k,\end{equation}
			where \begin{equation}\label{variance}s^2=\frac{1}{2}\sum\limits_{x\le p\le y} \left|a_p\right|^2.\end{equation}\end{prop}
		Note that the right-hand side of \eqref{bessapp} is the moment of a Gaussian random variable with mean $0$ and variance $s^2$.
		Substituting the values of $a_p$ from \eqref{achoice}, and setting $x=q_n, y=q_m$ in \eqref{bessapp}, an application of Mertens's estimate shows that:
		\begin{equation}s^2=\sum\limits_{q_m\le p\le q_n} \frac{1}{2p}=\frac{(\log \log q_m-\log \log q_n)}{2}+o(1).\end{equation}
		Hence, \begin{equation}
			\E\left[\left|T_1\right|^{2k}\right]\le \frac{(2k)!}{2^{2k} k!}(\log \log q_m-\log \log q_n)^k.
		\end{equation}
		We can use a similar method to bound $\E[\left|T_2\right|^{2k}]$.  We take $a(p)=\frac{1}{p\sqrt{2}}$ in  \eqref{bessapp}.
		Since $T_2^k$ has length $\le q^{\frac{1}{3}}$ we may pass to the random variables $X$.
		We obtain:
		\begin{equation}\E\left[\left|S_2\right|^{2k}\right]\le \frac{(2k)!}{2^{2k} k!} \left(\sum\limits_{p\in (q_n,q_m]} \frac{1}{p^2}\right)^k.\end{equation}
		
		Another application of H{\"o}lder's inequality shows: \begin{equation}\label{squareerror}\E\left[\left|S_{m}-S_n\right|^{2k}\right]\le \frac{(2k)!}{2^{2k} k!}(m-n)^k \left(1+O\left(\frac{\sum\limits_{p \in (q_n,q_m]}{\frac{1}{p^2}}}{\sum\limits_{p \in (q_n,q_m]}{\frac{1}{p}}}\right)\right)^k.\end{equation}
		We show the error term is negligible.
		If $n\ge 1$, then \begin{equation}\frac{\sum\limits_{p \in (q_n,q_m]}{\frac{1}{p^2}}}{\sum\limits_{p \in (q_n,q_m]}{\frac{1}{p}}}<\frac{1}{q_n}.\end{equation} 
		Since $k\le \log q,$ the error term in \eqref{squareerror} is $\ll 1+O\left(\frac{\log q}{q_n}\right)\ll 1$.
		If $n=0$, then $\sum\limits_{p \in (q_n,q_m]}{\frac{1}{p^2}}<1$, whilst \begin{equation}\sum\limits_{p \in (q_n,q_m]}{\frac{1}{p}}\gg \log \log q_m\gg \log\log q.\end{equation}
		Hence, we see $\left|\frac{S_2}{S_1}\right|\ll  \frac{1}{\log\log q}$.
		By the choice of $k$, we have $(q_m)^k \le q$, so that $k\le \log_3 q.$
		In either case, we obtain  \begin{equation}
			\left(1+O\left(\frac{\sum\limits_{p \in (q_n,q_m]}{\frac{1}{p^2}}}{\sum\limits_{p \in (q_n,q_m]}{\frac{1}{p}}}\right)\right)^k\ll 1,
		\end{equation}
		so the sums of squares of primes doesn't change the order of the sum and we conclude.\end{proof}
	
	\section{Twisted second moments in the {$q$}-aspect}\label{Twistsec}
	This section is dedicated to the proof of Theorem  \ref{realtwistthm}. We first need an expression for the mean value of $|M_1\dots M_l F|^2$ given in Proposition \ref{ourtwistthm}.
	
	Let $F$ be an $l$-sufficient function, as defined in \eqref{lsuff}. We use the same notation as in the proof of Lemma \ref{Realsplit}.
	By Equation \eqref{jointK}, we may write $F$ as \begin{equation}\label{defF}
		F=\mathlarger{\sum}\limits_{u,v} \frac{\tilde{C}_{u,v}}{\sqrt{uv}} \chi(u)\overline{\chi(v)}.
	\end{equation}
	By construction of the coefficients, we have  
	\begin{equation}
		\label{lenf}
		u,v\le\exp\left(\mathlarger{\sum}^l_{j=1} 10(n_j-n_{j-1})^{10^4}\log q_j\right)\ll  q^{10^{-3}}, 
	\end{equation}
	so that we may apply Lemma \ref{realpartexp} to each $F_j$ for $1\le j\le l$, and Lemma  \ref{Realsplit} to $F$.

	By Lemma \ref{Realsplit}, the second moment of $F$ may be expressed as:
	\begin{equation}\label{Fmom}
		\E_{\oplus}\left[F^2\right]=\mathlarger{\sum}\limits_{u,v} \frac{\left|\tilde{C}_{u,v}\right|^2}{uv}.
	\end{equation}
	We need to bound the twisted mollified moment by $\frac{\log q}{\log q_l}\E_{\oplus}\left[F^2\right]$.
	Let $\tilde{M}=M_1....M_l$
	and \begin{equation}\label{defM}
		M=\tilde{M}F
	\end{equation}
	be the mollified twist.
	Our restrictions on the coefficients of $M_l$ and $F$ ensure we may write \begin{equation}\label{Mcoef}M(\chi)=\sum^{\mathbf{M}}_{j,k=1} \frac{x_{j,k}}{\sqrt{jk}} \chi(j)\overline{\chi(k)},\end{equation} 
	for some choice of coefficients $x_{j,k}$ integers $j$ and $k$ up to a finite \textit{length} $\mathbf{M}. $ We define an array $X=(x_{j,k})_{j,k}.$
	\begin{prop}\label{ourtwistthm}
		Let $M$ and $X$ be defined as above. Then
		
		\begin{equation}\label{twistmomform}
			\E_{+}\left[\left|M\right|^2\right]\ll\varphi(q)q^{-1}\mathscr{Q}(X)+O\left(q^{-\frac{1}{3}} \max_{u,v}  \left|\tilde{C}_{u,v}\right|^2	\right),
		\end{equation}
	\end{prop}
	where \begin{equation}\label{Qformula}
		\mathscr{Q}(X)=\sum^{\mathbf{M}}_{\substack{j_1,j_2,k_1,k_2=1\\ (j_1,q)=(j_2,q)=(k_1,q)=(k_2,q)=1}}  x_{j_1,k_1}\overline{x_{j_2,k_2}}\frac{(j_1k_2,j_2k_1)}{j_1j_2k_1k_2}\log \left(\frac{R^2 (j_1k_2,j_2k_1)^2}{j_1j_2k_1k_2}\right),
	\end{equation}
	and  $R=R(q)$ is defined by: 
	\begin{equation}\label{Rdef}
		\log R=\frac{1}{2} \log \left(\frac{q}{\pi}\right)+\frac{1}{2}\psi\left(\frac{1}{4}\right)+\gamma+\eta(q).
	\end{equation}
	
	\begin{proof}
		In \cite{Sarnak}, (Section $5$, Equation $5.5$) they produce an estimate for the Fourier transform:
		\begin{equation} \label{3.6}
			B(m_1,m_2)=\mathlarger{\sum}^+_{\chi \text{(mod $q$)}} \chi(m_1)\overline{\chi(m_2)} \left|L\left(\chi,\frac{1}{2}\right)\right|^2. 	
		\end{equation}
		If $(c,q)=1$, then \begin{equation}\label{copfactor}
			B(cm_1,cm_2)=B(m_1,m_2),
		\end{equation}
		so in order to understand $B(m_1,m_2)$, it suffices to consider the case where $(m_1,m_2)=1$.
		In \cite{Sarnak}, they bound the  terms $B(m_1,m_2)$ in the following lemma. 
		\begin{lem}[Lemma 3.2 in \cite{Sarnak}] 
			Suppose $(m_1,m_2)=1$ and let $B(m_1,m_2)$ be as defined in \eqref{3.6}.
			Then
			\begin{equation}\label{Bsum}
				B(m_1,m_2)=\mathlarger{\sum}_{vw=q}\mu(v)\varphi(w)\mathlarger{\sum}^\star_{m_1n_1\equiv m_2n_2 \text{ mod }w} W(\pi n_1n_2/q),\end{equation}
			\textit{	where }\begin{equation}
				W(y)=\frac{1}{2\pi i}\int^{1+i\infty}_{1-i\infty} \Gamma^2 \left(\frac{s}{2}+\frac{1}{4}\right) G^2(s)s^{-1} y^{-s} ds,
			\end{equation}\textit{ and $G(s)$ is a holomorphic function in the region $\left|\Re(s)\right|<1$ such that \begin{equation}
					G(s)=G(-s), G\left(\frac{1}{2}\right)=G\left(-\frac{1}{2}\right)=0, G(0)\Gamma\left(\frac{1}{4}\right)=1, 
				\end{equation}
				and  }\begin{equation}
				G(s)\ll |s|^{-3}e^{\frac{\pi|s|}{4}}.
			\end{equation}
		\end{lem}
		In \eqref{Bsum}, the $\star$ restricts the summation to numbers coprime to $q$.
		Moreover, the dominant term is the diagonal contribution to $B(m_1,m_2)$ in \eqref{Bsum}  when $m_1n_1=m_2n_2$, which by their Lemma 4.1 is:
		\begin{equation}
			B_0(m_1,m_2)=\frac{\varphi^+(q)\varphi(q)}{q\sqrt{m_1m_2} }\log\left({\frac{R^2}{m_1m_2}}\right)+O(\tau(q)q^{\frac{1}{2}}).
		\end{equation}

		The error terms from the other terms contribute at most $\beta(m_1,m_2)$,	where 
		\begin{equation} \label{4.8}
			\beta(m_1,m_2)=\mathlarger{\sum}\limits_{\substack{m_1n_1\ne m_2n_2\\ (n_1,q)=(n_2,q)=1}} (m_1n_1\pm m_2n_2,q) (n_1n_2)^{-\frac{1}{2}} \left|W\left(\frac{\pi n_1n_2}{q}\right)\right|, 
		\end{equation}
		and the sums for the different signs for $\pm$ are taken separately.
		Using \eqref{4.8}, they obtain:\begin{equation}\label{FT}
			B(m_1,m_2)=\frac{\varphi^+(q)\varphi(q)}{q\sqrt{m_1m_2}}\log \left(\frac{R^2}{m_1m_2}\right)+O(\beta(m_1,m_2)+\tau(q)q^{\frac{1}{2}}), 
		\end{equation}
		where $R$ is as defined in \eqref{Rdef}.
		In order to get a sufficient error bound when we apply \eqref{FT} to the terms forming $M$, we require a bound on the values of $j$ and $k$ which contribute to the sum in 
		\eqref{Mcoef}.
		\begin{prop}
			If $\max\{j,k\}\gg q^{\frac{1}{100}}$, then $x_{j,k}=0$
		\end{prop} 
		\begin{proof} Non-zero coefficients in $M$ come from products of terms of $M_1...M_l$ and $F.$ By construction,  for each $1\le j\le l,M_j$ has length $\le \exp( (n_j-n_{j-1})^{10^5} \log q_j)$,
			and so by construction of the values $q_l$, $M_1...M_l$ has length $
			\ll  \exp\left((n_{l+1}-n_l)^{10^5}\right)e^{n_l},
			$
			which is $		\ll  q^{10^{-4}}$.
			
			By Equations \eqref{defF} and \eqref{lenf}, we may write 
			\begin{equation}F=\mathlarger{\sum}\limits_{
					u,v\ll q^{10^{-3}}} \frac{\tilde{C}_{u,v}}{\sqrt{uv}} \chi(u)\overline{\chi(v)}. 
			\end{equation}
			so that we may take \begin{equation}\label{Mbound}
				\mathbf{M}\ll  q^{\frac{1}{100}}.	
			\end{equation}
		\end{proof}

		For the context of $l$-sufficient functions, since the coefficients of $M_l$ have size at most $1$, and $M_1...M_l$ has length $\ll q^{10^{-4}}$,  the coefficients of $M$ satisfy:
		\begin{equation}\label{coefbound}
			\max_{j,k}\frac{\left|x_{j,k}\right|}{\sqrt{jk}}\ll q^{{10}^{-4}}\max_{u,v} \frac{\left|\tilde{C}_{u,v}\right|}{\sqrt{uv}}.
		\end{equation}
		
		Using \eqref{Mcoef}, we may write
		\begin{equation}
			\left|M\right|^2=\sum^{\mathbf{M}}_{\substack{j_1,j_2,k_1,k_2=1\\ (j_1,q)=(j_2,q)=(k_1,q)=(k_2,q)=1}} \frac{\chi(j_1k_2)\overline{\chi(j_2k_1)}}{\sqrt{j_1k_1j_2k_2}} x_{j_1,k_1}\overline{x_{j_2,k_2}}.
		\end{equation}
		The expression in \eqref{FT} for $B(m_1,m_2)$ is only valid when $m_1$ and $m_2$ are coprime. 
		To put all the non-zero terms in the above expression in this form, we rewrite it as:
		\begin{equation}
			\left|M\right|^2=\sum^{\mathbf{M}}_{\substack{j_1,j_2,k_1,k_2=1\\ (j_1,q)=(j_2,q)=(k_1,q)=(k_2,q)=1}} \chi\left(\frac{j_1k_2}{(j_1k_2,j_2k_1)}\right)\overline{\chi\left(\frac{j_2k_1}{(j_1k_2,j_2k_1)}\right)} \frac{x_{j_1,k_1}\overline{x_{j_2,k_2}}}{\sqrt{j_1k_1j_2k_2}}.
		\end{equation}
		Using \eqref{FT}, we can write:
		\begin{equation}
			\begin{aligned}\label{realsp}
				&\E_{+}\left[\left|LM\right|^2\right]\ll\sum^{\mathbf{M}}_{\substack{j_1,j_2,k_1,k_2=1\\ (j_1,q)=(j_2,q)=(k_1,q)=(k_2,q)=1}}  x_{j_1,k_1}\overline{x_{j_2,k_2}}
				\Big\{\frac{\varphi(q)(j_1k_2,j_2k_1)}{q\sqrt{j_1j_2k_1k_2}}\log \left(\frac{R^2 (j_1k_2,j_2k_1)^2}{j_1j_2k_1k_2}\right)\\
				&\hspace{4cm}+ O\Big(\frac{\beta\left(\frac{j_1k_2}{(j_1k_2,j_2k_1)},\frac{j_2k_1}{(j_1k_2,j_2k_1)}\right)+\tau(q)q^{\frac{1}{2}}}{\varphi(q)\sqrt{j_1k_1j_2k_2}}\Big)\Big\}
			\end{aligned}
		\end{equation}
		
		If we put $X=(x_{j,k})_{j,k}$,
		then we may rewrite the right-hand side of \eqref{realsp} as
		\begin{equation}\label{realerror}
			\varphi(q)q^{-1}\mathscr{Q}(X)+O\left(\frac{\max_{j,k} \frac{|x_{j,k}|^2}{jk}}{\varphi(q)}\left(\mathlarger{\sum}\limits_{j_1,j_2,k_1,k_2} \tau(q) q^{\frac{1}{2}}+\beta\left(\frac{j_1k_2}{(j_1k_2,j_2k_1)},\frac{j_2k_1}{(j_1k_2,j_2k_1)}\right)\right)\right).	
		\end{equation}

		Using the relation \eqref{coefbound}, we may re-express the error bound as
		\begin{equation}	\label{realmid}
			\ll O\left(\frac{q^{10^{-3}}\max_{u,v} |\tilde{C}_{u,v}|^2}{\varphi(q)}\left(\sum^{\mathbf{M}}_{\substack{j_1,j_2,k_1,k_2=1\\ (j_1,q)=(j_2,q)=(k_1,q)=(k_2,q)=1}} \tau(q) q^{\frac{1}{2}}+\beta\left(\frac{j_1k_2}{(j_1k_2,j_2k_1)},\frac{j_2k_1}{(j_1k_2,j_2k_1)}\right)	\right)
			\right).
		\end{equation}

		We need to bound \begin{equation}
			\sum^{\mathbf{M}}_{\substack{j_1,j_2,k_1,k_2=1\\ (j_1,q)=(j_2,q)=(k_1,q)=(k_2,q)=1}}\beta\left(\frac{j_1k_2}{(j_1k_2,j_2k_1)},\frac{j_2k_1}{(j_1k_2,j_2k_1)}\right).
		\end{equation}
		We recall from \eqref{4.8}, if $(m_1,m_2)=1$ then,
		\begin{equation} 
			\beta(m_1,m_2)=\mathlarger{\sum}\limits_{\substack{m_1n_1\ne m_2n_2\\ (n_1,q)=(n_2,q)=1}} (m_1n_1\pm m_2n_2,q) (n_1n_2)^{-\frac{1}{2}} \left|W\left(\frac{\pi n_1n_2}{q}\right)\right|. 
		\end{equation}
		Hence, \begin{align}\label{betasum}
			&\sum^{\mathbf{M}}_{\substack{j_1,j_2,k_1,k_2=1\\ (j_1,q)=(j_2,q)=(k_1,q)=(k_2,q)=1}}\beta\left(\frac{j_1k_2}{(j_1k_2,j_2k_1)},\frac{j_2k_1}{(j_1k_2,j_2k_1)}\right)\\\nonumber&=\mathlarger{\sum}\limits_{(n_1,q)=(n_2,q)=1}(n_1n_2)^{-\frac{1}{2}} \left|W\left(\frac{\pi n_1n_2}{q}\right)\right|\sum^{\mathbf{M}}_{\substack{j_1,j_2,k_1,k_2=1\\ n_1j_1k_2\ne n_2j_2k_1}}\left(\frac{n_1j_1k_2 \pm n_2j_2k_1}{(j_1k_2,j_2k_1)},q\right).\end{align}
		Since $j_1,j_2, k_1,k_2n_1$ and $n_2$ are coprime to $q$, the inner sum:
		\begin{equation}
			\sum^{\mathbf{M}}_{\substack{j_1,j_2,k_1,k_2=1\\ n_1j_1k_2\ne n_2j_2k_1}}\left(\frac{n_1j_1k_2 \pm n_2j_2k_1}{(j_1k_2,j_2k_1)},q\right)=\sum^{\mathbf{M}}_{\substack{j_1,j_2,k_1,k_2=1\\ n_1j_1k_2\ne n_2j_2k_1}}(n_1j_1k_2
			\pm n_2j_2k_1,q)\end{equation}
		may be bounded as $\ll \mathlarger{\sum}\limits_{d|q} \frac{\mathbf{M}^4}{d}\ll \mathbf{M}^4 \tau(q)$. 
		For Equation 4.8 in \cite{Sarnak}, the authors use $W(y)\ll (1+y)^{-1}$.
		Hence the sum in \eqref{betasum} may be bounded by
		\begin{equation}
			\begin{aligned}
				\mathbf{M}^4\tau(q)	\mathlarger{\sum}\limits_{n_1,n_2}(n_1n_2)^{-\frac{1}{2}}\left(1+\frac{n_1n_2}{q}\right)^{-1}
				&\ll \mathbf{M}^4 \left(\tau(q)\mathlarger{\sum}\limits_{\substack{n_1,n_2\\ n_1n_2<q}}(n_1n_2)^{-\frac{1}{2}}+q\mathlarger{\sum}\limits_{n_1n_2>q}(n_1n_2)^{-\frac{3}{2}}\right) \\
				&\ll \mathbf{M}^4 \tau(q) \sqrt{q} \log q.
			\end{aligned}
		\end{equation}
		
		Using this bound for the error term in \eqref{realmid} for the expression in \eqref{realerror}, the bound in \eqref{Mbound}, and the bounds \begin{equation}
			\frac{q}{\varphi(q)}, \tau(q)\ll_\epsilon q^\epsilon,
		\end{equation} we obtain 
		\begin{equation}\label{realend}
			\E_{+}\left[\left|M\right|^2\right]\ll\varphi(q)q^{-1}\mathscr{Q}(X)+O\left(q^{-\frac{1}{3}} \max_{u,v} |\tilde{C}_{u,v}|^2\right).
		\end{equation}
		This completes the proof of Proposition  \ref{ourtwistthm}.
	\end{proof}
	
	\begin{proof}[Proof of Theorem \ref{realtwistthm}]
		Having proven Proposition \ref{ourtwistthm}, we proceed with the proof of Theorem \ref{realtwistthm}.	We want to calculate $\mathscr{Q}(X)$ for the value of $X$ determined by the coefficients of $M$ as defined in \eqref{defM}.
		
		We recall from Equation \eqref{Fmom}, the following formula for the expectation of the function $F$: 
		\begin{equation}\label{Fcopy}
			\E_{\oplus}\left[F^2\right]=\mathlarger{\sum}\limits_{\substack{p|uv\implies p\le q_l\\
					(u,v)=1}}\frac{ |\tilde{C}_{u,v}|^2}{uv}.
		\end{equation}
		
		Combining Equations \eqref{coefbound} and \eqref{Fcopy} with the bounds \begin{equation}
			\frac{q}{\varphi(q)}, \frac{q}{\varphi^+(q)}, \tau(q)\ll _\epsilon q^\epsilon,
		\end{equation}
		for any $\epsilon>0$, we see that the error term in \eqref{twistmomform} is \begin{equation}\label{Qerror}
			o\left(\E_{\oplus}\left[\left|F\right|^2\right]\right).
		\end{equation}
		We will show that \begin{equation}\label{Qbound}
			\mathscr{Q}(X)\ll \frac{\log q}{\log q_l} \E_{+}\left[\left|F\right|^2\right],
		\end{equation}
		which upon substitution into \eqref{twistmomform} will complete the proof of Theorem \ref{realtwistthm}. 
		
		We need to substitute the choice of $M$ into \eqref{twistmomform}. First, we rewrite the expression for $\mathscr{Q}(X)$ in \eqref{Qformula} to group terms with coming from the same term in $Q(\chi)$ together.
		Let $\tilde{M}$ have  coefficients given by:
		\begin{equation}
			\tilde{M}(\chi)=\mathlarger{\sum}\limits_{f} \frac{e_f\chi(f)}{m^{\frac{1}{2}}}.
		\end{equation}
		Expanding $M=\tilde{M}F$, we see
		\begin{equation}
			x_{j,k}=\mathlarger{\sum}\limits_{uf=j}\tilde{C}_{u,k}e_f.
		\end{equation}
		Substituting this value into \eqref{Qformula} yields:
		\begin{equation}\label{Qunsplit}
			\mathscr{Q}(X)=\mathlarger{\sum}_{u_1,k_1,u_2,k_2} \tilde{C}_{u_1,k_1}\overline{\tilde{C}_{u_2,k_2}}\quad \mathlarger{\sum}_{f_1,f_2}\frac{ e_{f_1}{e_{f_2}}}{[u_1k_2f_1,u_2k_1f_2]} \log\left(\frac{R^2(u_1k_2f_1,u_2k_1f_2)^2}{u_1k_2f_1u_2k_1f_2}\right).
		\end{equation}
		
		For fixed choices of $u_1,k_1,u_2$ and $k_2$, we can separate the log factor in \eqref{Qunsplit} into two terms:  \begin{equation}\log\left(\frac{R^2(u_1k_1f_2,u_2k_2f_1)^2}{u_1k_1f_2u_2k_2f_1}\right)=\log\left(\frac{R^2(u_1k_2,u_2k_1)^2}{u_1k_1u_2k_2}\right)-\log\left(\frac{f_1f_2(u_1k_2,u_2k_1)^2}{(u_1k_2f_1,u_2k_1f_2)^2}\right).
		\end{equation}
		If we take this first term on the right hand side, which is constant in $f_1$ and $f_2$, out, then the remaining factor inside the logarithm is doubly multiplicative in $f_1$ and $f_2$ for a fixed choice of $c_1$ and $c_2$, which makes it easier to sum over.
		We write \eqref{Qunsplit} as:
		\begin{equation}\label{QSsplit}
			\mathscr{Q}(X)=P_1-P_2,
		\end{equation}
		where \begin{equation}\label{defS1}
			P_1=\mathlarger{\sum}_{u_1,k_1,u_2,k_2} \tilde{C}_{u_1,k_1}\overline{\tilde{C}_{u_2,k_2}}\log\left(\frac{R^2(u_1k_2,u_2k_1)^2}{u_1k_1u_2k_2}\right)\mathlarger{\sum}_{f_1,f_2}\frac{ e_{f_1}{e_{f_2}}}{[u_1k_2f_1,u_2k_1f_2]},
		\end{equation}
		and 
		\begin{equation}\label{S2def} P_2=\mathlarger{\sum}_{u_1,k_1,u_2,k_2 } \tilde{C}_{u_1,k_1}\overline{\tilde{C}_{u_2,k_2}}\mathlarger{\sum}_{f_1,f_2}\frac{ e_{f_1}{e_{f_2}}}{[u_1k_2f_1,u_2k_1f_2]} \log\left(\frac{f_1f_2(u_1k_2,u_2k_1)^2}{(u_1k_2f_1,u_2k_1f_2)^2}\right).
		\end{equation}
		
		We think of $P_1$ as being the sum where we take an approximate value for the logarithm factor by setting $f_1=f_2=1$, and $P_2$ being the oscillatory term which takes into account the variation in the logarithmic factor for different values of $f_1$ and $f_2$. 
		We bound both $P_1$ and $P_2$, and show: \begin{equation}
			P_1,P_2\ll \frac{\log q}{\log q_l}\E_{\oplus} \left[\left|F\right|^2\right].
		\end{equation}
		We first bound $P_1$ by taking the maximum value of the logarithm factor.
		From the definition of $R$, we see the logarithmic factor be bounded by $\log q$ and so 
		\begin{equation}
			\left|P_1\right|\le \mathlarger{\sum}_{u_1,k_1,u_2,k_2} \left|\tilde{C}_{u_1,k_1}\overline{\tilde{C}_{u_2,k_2}}\log\left(\frac{R^2(u_1k_2,u_2k_1)^2}{u_1k_1u_2k_2}\right)\right|\left|\mathlarger{\sum}_{f_1,f_2}\frac{ e_{f_1}{e_{f_2}}}{[u_1k_2f_1,u_2k_1f_2]}\right|
			\ll P_3,\end{equation}
		where 
		\begin{equation}\label{P3def}
			P_3=\log q  \mathlarger{\sum}_{u_1,k_1,u_2,k_2} \left|\tilde{C}_{u_1,k_1}\overline{\tilde{C}_{u_2,k_2}}\right|\left|\mathlarger{\sum}_{f_1,f_2}\frac{ e_{f_1}{e_{f_2}}}{[u_1k_2f_1,u_2k_1f_2]}\right|.
		\end{equation}
		We show that $P_3$ may be bounded by $\frac{\log q}{\log q_l}\E_{\oplus}[\left|Q\right|^2]$ in Section \ref{sect: prop S3lem} below.
		\begin{prop}\label{S3lem}
			Let $P_3$ be defined as above. Then \begin{equation}
				P_3\ll   \frac{\log q}{\log q_l}\E_{\oplus}[F^2].
			\end{equation}
		\end{prop}
		
		For $P_2$, we also want to bound the sum absolutely to begin with.
		We have from \eqref{S2def},
		\begin{equation}\label{S2new}\left|P_2\right|\le\mathlarger{\sum}_{u_1,k_1,u_2,k_2} \left|\tilde{C}_{u_1,k_1}\tilde{C}_{u_2,k_2}\right|\left|\mathlarger{\sum}_{f_1,f_2}\frac{ e_{f_1}{e_{f_2}}}{[u_1k_2f_1,u_2k_1f_2]} \log\left(\frac{f_1f_2(u_1k_2,u_2k_1)^2}{(u_1k_2f_1,u_2k_1f_2)^2}\right)\right|.
		\end{equation}
		
		We want to approximate the logarithm by multiplicative functions, to make them easier to sum.
		Observe that,  given $T>1$, as $\beta\to 0$ we have by Taylor's Theorem,
		\begin{equation}\label{loga}
			\log t=\frac{t^{i\beta}-t^{-i\beta}}{2i\beta}+O(\beta),\end{equation}
		where the error decays uniformly in $\beta$, for all $1\le t\le T.$
		This will enable us to bound $P_2$.
		
		We substitute $t$ for  \begin{equation} \frac{f_1f_2(u_1k_2,u_2k_1)^2}{(u_1k_2f_1,u_2k_1f_2)^2}	
		\end{equation} in the definition of $P_2$ in \eqref{S2def}. We see $t\ge 1$, whilst the restrictions on prime factorisations ensure that if $e_{f_1}e_{f_2}\ne 0$, then $t$ can be bounded by some $T$ depending only on $q$. Hence, $P_2$ may be approximated by taking the limit as $\beta\to 0$ in \eqref{loga}. Indeed, if we set:   
		\begin{align}\label{P4def}
			P_4(\beta)=\frac{1}{2\beta}&\mathlarger{\sum}_{u_1,k_1,u_2,k_2} \left|\tilde{C}_{u_1,k_1}\overline{\tilde{C}_{u_2,k_2}}\right|\\\nonumber&\left|\mathlarger{\sum}_{f_1,f_2}\frac{ e_{f_1}{e_{f_2}}}{[u_1k_2f_1,u_2k_1f_2]} \left(\left(\frac{f_1f_2(u_1k_2,u_2k_1)^2}{(u_1k_2f_1,u_2k_1f_2)^2}\right)^{i\beta}-\left(\frac{f_1f_2(u_1k_2,u_2k_1)^2}{(u_1k_2f_1,u_2k_1f_2)^2}\right)^{-i\beta}\right)\right|, 
		\end{align}
		then the bound in \eqref{S2new} becomes $\left|P_2\right|\le\lim_{\beta\to 0^+} P_4(\beta)$.
		We show in Section \ref{sect: S4lem}:
		\begin{prop}\label{S4lem} 	Let $P_4$ be defined as above. Then
			\begin{equation}\lim_{\beta\to 0^+}P_4(\beta)\ll  \frac{\log q}{\log q_l}\E_{\oplus}\left[F^2\right].\end{equation} 
		\end{prop}
		Combining Propositions  \ref{S3lem} and \ref{S4lem} completes the proof of Theorem \ref{realtwistthm}.
	\end{proof}
	
	\subsection{Proof of Proposition \ref{S3lem}}
	\label{sect: prop S3lem}
	The formula for $P_3$ is used for estimating the sum $P_1$ in \eqref{QSsplit}. The bound $P_3\ge P_1$ was obtained by bounding  the absolute value of the logarithmic factor by its maximum, which is smaller than $\log q$. The expression for $P_3$ splits into products over contributions of expectations of the components of the Dirichlet polynomial supported on  the integers with primes in the interval $(q_{i-1},q_i]$, for $1\le i\le l$.
	Since $F$ is $l$-sufficient, the contributions of the different intervals are weakly dependent, and the expectation of the $P_3$ may be bounded by the product of the expectations from each interval.
	Indeed, we have:
	\begin{equation} \label{S3prod}P_3=\log q\prod_i N_i, \end{equation}
	where $N_i$ is the contribution to $P_3$ from the $F_i$ component in the factorisation of $F=\prod^l_{i=1}F_i$:
	\begin{equation}\label{S3rest}
		N_i= \Big|\mathlarger{\sum}_{\substack{p|f_1,f_2\implies p\in (q_{i-1},q_i]\\
				\Omega_i(f_1),\Omega_i(f_2)\le (n_i-n_{i-1})^{10^5}}} \frac{ \mu(f_1)\mu(f_2)}{[u_1k_2f_1,u_2k_1f_2]}\Big| \mathlarger{\sum}_{u_1,k_1,u_2,k_2} \left|C_{u_1,k_1, B^{(i)}, K_i}C_{u_2,k_2,B^{(i)}, K_i}\right|.\end{equation}
	In Lemma \ref{S3approx} below, we show we may lift the restriction on the number of prime factors of coefficients $f_1$ and $f_2$ of the mollifier $M_l$ with negligible error; in Lemma \ref{S3eval} we bound the sum with the restriction removed.
	
	We use Rankin's trick to show that we can remove the restriction on the number of prime factors of $f_1$ and $f_2$, the terms coming from the mollifier $\tilde{M},$ with negligible error. If we were able to remove this restriction completely, then $\tilde{M}$ would completely mollify the component of the Euler product of the Dirichlet L-function coming from primes in the interval $(q_{i-1},q_{i}]$; this next lemma essentially shows that the factors of mollifier $M_1,...,M_l$ are long enough to have this mollification effect. Lemma \ref{S3eval} then evaluates the sum with this restriction removed.
	\begin{lem}\label{S3approx}
		For $0\le i\le l$, the expression in \eqref{S3rest} is 
		\begin{align}\label{S3unrest}
			\mathlarger{\sum}_{u_1,k_1,u_2,k_2} \left|C_{u_1,k_1, B^{(i)}, K_i}C_{u_2,k_2,B^{(i)}, K_i}\right|&\left|\mathlarger{\sum}_{p|f_1,f_2\implies p\in (q_{i-1},q_i]
			} \frac{ \mu(f_1)\mu(f_2)}{[u_1k_2f_1,u_2k_1f_2]}\right|\\\nonumber&+O\left(e^{-100(n_i-n_{i-1})}\E_{\oplus}\left[\left|F_i\right|^2\right]\right).
		\end{align}
	\end{lem}
	\begin{proof}	
		The contribution of those values of $f_1$ where $\Omega_i(f_1)> (n_i-n_{i-1})^{10^5},$ which are excluded by the restriction on the length of the mollifier, to \eqref{S3unrest} may be bounded for any $\rho>0$ as:
		\begin{equation}\label{CHer}
			\begin{aligned}
				&\le e^{-\rho(n_i-n_{i-1})^{10^5}}\times\\
				& \mathlarger{\sum}_{u_1,k_1,u_2,k_2} \left|C_{u_1,k_1, B^{(i)}, K_i}C_{u_2,k_2,B^{(i)}, K_i}\right|\mathlarger{\sum}_{\substack{p|f_1\implies p\in (q_{i-1},q_i]\\p|f_2\implies p\in (q_{i-1},q_i]
				}} e^{\rho\Omega(f_1)}\frac{\mu^2(f_1)\mu^2(f_2)}{[u_1k_2f_1,u_2k_1f_2]}.
			\end{aligned}
		\end{equation}
		We sum over the possible indices of each prime $p$ in the interval in the prime factorisations of $f_1$ and $f_2$, which both must be squarefree to contribute to the sum. Say the exponent is $d_{p,i}=v_p(f_i)\in \{0,1\}$.
		Then, for fixed values of $u_1,u_2,k_1$ and $k_2$, the inner sum over $f_1$ and $f_2$ in \eqref{CHer} is:
		\begin{equation}\label{termcoef}\prod_{p\in(q_{i-1},q_i]} \mathlarger{\sum}_{d_{p,1},d_{p,2}\in \{0,1\}} \frac{e^{\rho \omega_i(p^{d_{p,1})}}}{p^{\max\{(d_{p,1}+v_p(u_1k_2)), d_{p,2}+v_p(u_2k_1)\}}}\end{equation}
		For a given prime $p$ in the interval $(q_{i-1},q_i],$ the term in the inner product may be expressed as:
		\begin{equation}\frac{1}{p^{\max\{v_p(u_1k_2), v_p(u_2k_1)\}}} \mathlarger{\sum}_{d_{p,1},d_{p,2}\in \{0,1\}} \frac{e^{\rho \omega_i(p^{d_{p,1})}}}{p^{\max\{d_{p,1}+v_p(u_1k_2), d_{p,2}+v_p(u_2k_1)\}-\max\{v_p(u_1k_2), v_p(u_2k_1)\}}}.\end{equation}
		
		Considering the possible values for  the exponent $$\max\{d_{p,1}+v_p(u_1k_2)), d_{p,2}+v_p(u_2k_1)\}-\max\{v_p(u_1k_2), v_p(u_2k_1)\},$$this term associated to the prime $p$ may be bounded by:
		\begin{equation}\label{pbound}
			\frac{1}{{p^{\max\{v_p(u_1k_2), v_p(u_2k_1)\}}}}\times \begin{cases} 1+\frac{100 e^\rho}{p} & v_p(u_1k_2)= v_p(u_2k_1)\\
				100 e^\rho & v_p(u_1k_2)\ne v_p(u_2k_1)\end{cases}.
		\end{equation}
		Since $u_1, u_2,k_1$ and $k_2$ have relatively few prime factors, we expect for typical primes in the interval, $v_p(u_1k_2)=v_p(u_2k_1)=0$, so that  \eqref{pbound} only yields the larger bound $100 e^\rho$ on a few primes.
		Note that the contribution of $p\in (q_{i-1},q_{i}]$ with $v_p(u_1k_2)= v_p(u_2k_1)$ may be bounded by taking the product of the bound in \eqref{pbound} over all primes in the interval.
		Thus, their contribution is at most
		
		\begin{equation}
			\prod_{ p\in (q_{i-1},q_i]} 1+\frac{100 e^\rho}{p}.
		\end{equation}
		An application of the Prime Number Theorem shows that the contribution of these typical primes may be bounded by:
		\begin{equation}
			\left(\frac{\log q_i}{\log q_{i-1}}\right)^{100e^\rho}=e^{100e^\rho(n_i-n_{i-1})}	.
		\end{equation}

		Hence, if we define a function $\alpha$ to collate the contribution of all the bounds in \eqref{pbound} used for each prime in the interval, then the coefficient for each term in the polynomial, which we are calculating in \eqref{termcoef} may be bounded by: 
		\begin{equation}\label{alphacoef}
			\frac{\alpha(u_1k_2,u_2k_1)}{[u_1k_2,u_2k_1]} e^{100e^\rho(n_i-n_{i-1})},
		\end{equation}
		where $\alpha$ is a doubly multiplicative function such that:
		\begin{equation}\label{alpha}
			\alpha(p^n,p^n)=1+\frac{100 e^\rho}{p}\text{ for all $n\ge 1$, and } \alpha(p^n,p^m)=100e^\rho \text{ for $n,m\ge 0$ with $n\ne m$.}
		\end{equation}
		
		Returning to \eqref{CHer}, the effect of adding in the extra values of $f_1$ with many prime factors may be bounded by
		\begin{equation}\label{preCS}
			e^{-\rho(n_i-n_{i-1})^{10^5}+100e^\rho(n_i-n_{i-1})} \mathlarger{\sum}_{u_1,k_1,u_2,k_2} \left|C_{u_1,k_1, B^{(i)}, K_i}C_{u_2,k_2,B^{(i)}, K_i}\right| \frac{\alpha(u_1k_2,u_2k_1)}{[u_1k_2,u_2k_1]}.
		\end{equation} 
		Using the Cauchy-Schwarz inequality:
		\begin{equation}
			\left|C_{u_1,k_1, B^{(i)}, K_i}C_{u_2,k_2,B^{(i)}, K_i}\right|\le \frac{1}{2}\left(\left|C_{u_1,k_1, B^{(i)}, K_i}\right|^2+\left|C_{u_2,k_2,B^{(i)}, K_i}\right|^2\right),
		\end{equation} and dropping the restriction on $\Omega_i(u_1)$ and $\Omega_i(k_2)$,
		we can  bound the sum in  \eqref{preCS} by:
		
		\begin{equation}\label{tosum}
			\mathlarger{\sum}_{u_2,k_2} \left|C_{u_2,k_2,B^{(i)}, K_i}\right|^2\mathlarger{\sum}\limits_{\substack{p|u_1k_1\implies p\in (q_{i-1},q_i]\\
					(u_1,k_1)=(u_2,k_2)=1}} \frac{\alpha(u_1k_2,u_2k_1)}{[u_1k_2,u_2k_1]}.
		\end{equation}
		We seek to bound the inner sum, \begin{equation}\label{innertosum}
			\mathlarger{\sum}\limits_{\substack{p|u_1k_2\implies p\in (q_{i-1},q_i]\\
					(u_1,k_1)=(u_2,k_2)=1}} \frac{\alpha(u_1k_2,u_2k_1)}{[u_1k_2,u_2k_1]}.\end{equation}
		The restriction $(u_1, k_1)=(u_2,k_2)=1$ ensures that \begin{equation}[u_1k_2,u_2k_1]=[u_1,u_2][k_1,k_2].\end{equation}
		Meanwhile, using that for any prime $p\in(q_{i-1},q_i]$, for $v_p(u_1k_2)\ne v_p(u_2k_1)$ we must have either $v_p(u_1)\ne v_p(u_2)$ or $v_p(k_1)\ne v_p(k_2)$, we see that \begin{equation}\alpha(u_1k_2,u_2k_1)\le \alpha(u_1,u_2)\alpha(k_2,k_1).\end{equation}
		Hence we may bound the sum in \eqref{innertosum} for any values of $u_2$ and $k_2$ by:
		\begin{equation}\label{innerto2}
			\mathlarger{\sum}\limits_{p|u_1k_2\implies p\in (q_{i-1},q_i]} \frac{\alpha(u_1,u_2)\alpha(k_1,k_2)}{[u_1,u_2][k_1,k_2]}=\left(\mathlarger{\sum}\limits_{p|u_1\implies p\in (q_{i-1},q_i]} \frac{\alpha(u_1,u_2)}{[u_1,u_2]}\right)\left(\mathlarger{\sum}\limits_{p|k_2\implies p\in (q_{i-1},q_i]} \frac{\alpha(k_1,k_2)}{[k_1,k_2]}\right).
		\end{equation}
		
		We follow \cite{ABR20}, Section $8.2$, where we change their notation to write $\alpha$ for their function $f$. Using their bounds, we have for any integer $c_2$ whose prime factors all lie in the interval $(q_{i-1},q_i]:$
		\begin{equation}\label{sumbound}
			\mathlarger{\sum}_{p|c_1\implies p\in (q_{i-1},q_i]} \frac{\alpha(c_1,c_2)}{[c_1,c_2]}\le \frac{1}{c_2}\exp(200e^\rho(n_i-n_{i-1})+10^4 \rho(n_i-n_{i-1})^{10^4+1}.
		\end{equation}
		Combining \eqref{alphacoef}, \eqref{tosum} and \eqref{sumbound} we see that \eqref{CHer} may be bounded by:
		\begin{equation}\label{Cherprep}
			e^{-\rho(n_i-n_{i-1})^{10^5}+500e^\rho(n_i-n_{i-1})+2.10^4\rho(n_i-n_{i-1})^{10^4+1}}\mathlarger{\sum}_{u_2,k_2} \left|C_{u_2,k_2,B^{(i)}, K_i}\right|^2.
		\end{equation}
		For $1\le i\le l$ we have by \eqref{defn}, \begin{equation}
			n_i-n_{i-1}=\mathbf{s} (\log_{i+1}(q)-\log_{i+2}(q)),
		\end{equation}
		and so, by choice of $\mathbf{s}$, we have $n_i-n_{i-1}\ge 10^5$.
		Hence, picking $\rho=1000$ in \eqref{Cherprep}, we may bound \eqref{CHer} by
		$$e^{-100(n_i-n_{i-1})}\E_{\oplus}\left[|F_i|^2\right].$$ 
		
		A similar calculation shows we may lift the restriction on $\Omega(f_2)$, incurring an error of at most $e^{-100(n_i-n_{i-1})}\E_{\oplus}\left[|F_i|^2\right].$ This concludes the proof of Lemma \ref{S3approx}.
	\end{proof}
	It remains to estimate the main  term in the product for $P_3$ in \eqref{S3unrest}, i.e.,
	\begin{equation}\label{unrestmain}\mathlarger{\sum}_{\substack{u_1,k_1,u_2,k_2}} \left|C_{u_1,k_1, B^{(i)}, K_i}C_{u_2,k_2,B^{(i)}, K_i}\right|\left|\mathlarger{\sum}_{p|f_1,f_2\implies p\in (q_{i-1},q_i]
		} \frac{ \mu(f_1)\mu(f_2)}{[u_1k_2f_1,u_2k_1f_2]}\right|.\end{equation}
	\begin{lem}\label{S3eval}
		The value of \eqref{unrestmain} is \begin{equation}
			\E_{\oplus}\left[\left|F_i\right|^2\right]\times \prod_{p\in(q_{i-1},q_i]} (1-\frac{1}{p}).  
		\end{equation}
	\end{lem}
	\begin{proof}
		We consider the contribution to the sum for different values of $u_1,u_2,k_1$ and $k_2$. We recall that the coprimality conditions for the support of the summand mean that $u_1k_2=u_2k_1\iff u_1=u_2,k_1=k_2.$ Suppose $u_1k_2\ne u_2k_1$. We show the sum vanishes, so that only the diagonal terms from taking $u_1k_2=u_2k_1$ remain.
		Since $u_1k_2$ and $u_2k_1$ have distinct prime factorisations, there must be some prime which divides them to a different exponent. 
		Hence, \begin{equation}\exists p^\ast \in (q_{i-1},q_i]: v_{p^\ast}(u_2k_1)\ne v_{p^\ast}(u_1k_2).\end{equation}
		Without loss of generality, we may assume \begin{equation}\label{porder}v_{p^\ast}(u_2k_1)>v_{p^\ast}(u_1k_2).\end{equation}
		We consider the prime decomposition of $f_1,f_2,u_2k_1$ and $u_1k_2$ into powers of $p^\ast$ and a remainder.
		Write 
		$$
		f_1=f_1'(p^\ast)^{d_1}\quad
		f_2=f_2'(p^\ast)^{d_2}\quad
		u_2k_1=c_1(p^\ast)^{d_3}\quad
		u_1k_2=c_2(p^\ast)^{d_4},
		$$
		where $p^\ast \nmid f_1',f_2', c_1, c_2$ and $d_1,d_2\in \{0,1\}.$
		The factor multiplying $\left|C_{u_1,k_1, B^{(i)}, K_i}C_{u_2,k_2,B^{(i)}, K_i}\right|$ in \eqref{unrestmain} is:
		\begin{equation}\mathlarger{\sum}_{p|f_1',f_2'\implies p\in (q_{i-1},q_i]\setminus p^\ast} \frac{\mu(f_1')\mu(f_2')}{[c_1f_1', c_2f_2']}\mathlarger{\sum}^1_{d_1,d_2=0} \frac{(-1)^{d_1}(-1)^{d_2}}{(p^\ast)^{\max\{d_1+d_3, d_2+d_4\}}}.\end{equation}
		From \eqref{porder}, we see $d_3>d_4$, so that the inner sum over $(d_1,d_2)$ is:	
		\begin{equation}\label{vanish}
			\mathlarger{\sum}^1_{d_1=0}\frac{(-1)^{d_3}}{(p^\ast)^{d_1+d_3}}\mathlarger{\sum}^1_{d_2=0} (-1)^{d_2}=0,
		\end{equation}
		so the whole sum vanishes in the case $u_2k_1\ne u_1k_2.$ We are left with the sum over the diagonal terms $u_2k_1=u_1k_2$ in \eqref{unrestmain}, which is:
		\begin{equation}\mathlarger{\sum}_{\substack{p|u_2k_1\implies p\in (q_{i-1},q_i]\\ \Omega_i(u_2),  \Omega_i(k_1)\le 10(q_i-q_{i-1})^{10^4}}}\frac{\left|\gamma(u_2k_1)\right|^2}{u_2k_1}\left|\mathlarger{\sum}_{p|f_1,f_2\implies p\in (q_{i-1},q_i]} \frac{\mu(f_1)\mu(f_2)}{[f_1,f_2]} \right|.\end{equation}
		Since the functions are multiplicative, the factor coming from the sum over $f_1$ and $f_2$ is $\prod_{p\in(q_{i-1},q_i]} (1-\frac{1}{p})$,
		which multiplies each $u_2k_1$ term, completing the proof of Lemma \ref{S3eval} .
	\end{proof}
	
	Combining Lemma \ref{S3approx} and Lemma \ref{S3eval}, we see the contribution to the product in \eqref{S3rest} is:
	\begin{equation}
		N_i=\E_{\oplus}\left[\left|F_i\right|^2\right] \left(\prod_{p\in(q_{i-1},q_i]}(1-\frac{1}{p})+O(e^{-100(n_i-n_{i-1})}) \right).
	\end{equation}
	Using the Prime Number Theorem to estimate the product over $p$, this yields:
	\begin{equation}
		N_i=\E_{\oplus}\left[\left|F_i\right|^2\right] \prod_{p\in(q_{i-1},q_i]}(1-\frac{1}{p})(1+O(e^{-99(n_i-n_{i-1})})) .  
	\end{equation}
	Substituting this value of $N_i$ into \eqref{S3prod}, we take the product to get
	\begin{equation}P_3\ll\frac{\log q}{\log q_l}\prod^l_{i=1} \E_{\oplus}\left[\left|F_i\right|^2\right]=\frac{\log q}{\log q_l}\E_{\oplus}\left[\left|F\right|^2\right],\end{equation}
	where we used Lemma \ref{Realsplit} to multiply the expectations. This concludes the proof of Proposition \ref{S3lem}.
	
	\subsection{Proof of Proposition \ref{S4lem}}
	\label{sect: S4lem}
	We want to write $P_4(\beta)$ in terms of the sums of multiplicative functions, which can then be calculated accurately. As in the calculation of $P_3$, we first want to remove the restriction on the number of prime factors of $f_1$ and $f_2$, using Rankin's trick, and then evaluate the sum with these restrictions removed.
	
	Given integers $u_1,k_1,u_2,k_2$ with prime factors all at most $q_l$ such that $(u_1,k_1)=(u_2,k_2)=1$, define \begin{equation}\label{Hdef}
		H(\beta,u_1,k_1,u_2,k_2)= \mathlarger{\sum}_{\substack{p|f_1f_2\implies p\in (q_{m-1},q_m]}}\frac{ e_{f_1}{e_{f_2}}}{[u_1k_2f_1,u_2k_1f_2]} \left(\frac{f_1f_2(u_1k_2,u_2k_1)^2}{(u_1k_2f_1,u_2k_1f_2)^2}\right)^{i\beta},
	\end{equation}
	so that
	\begin{equation}	
		P_4(\beta)=\frac{1}{2\beta}\mathlarger{\sum}_{u_1,k_1,u_2,k_2} \left|\tilde{C}_{u_1,k_1}\tilde{C}_{u_2,k_2}\right|\left|H(\beta,u_1,k_1,u_2,k_2)-H(-\beta,u_1,k_1,u_2,k_2)\right|. 
	\end{equation}
	
	The expression for $H(\beta,u_1,k_1,u_2,k_2)$ in \eqref{Hdef} splits into the products of contribution supported on integers with prime factors all in the intervals  $(q_m,q_{m-1}]$ for each $1\le m\le l$. 
	Indeed, we may express Equation \eqref{Hdef} as:
	\begin{equation}
		H(\beta,u_1,k_1,u_2,k_2)=\prod^l_{m=1}R(m,\beta,u_1,k_1,u_2,k_2),
	\end{equation}
	where $R(m,\beta,u_1,k_1,u_2,k_2)$ is defined as:
	\begin{equation}\label{rmb}
		\mathlarger{\sum}_{\substack{p|f_1f_2\implies p\in (q_{m-1},q_m]}}\frac{ e_{f_1}{e_{f_2}}}{[u_1k_2f_1,u_2k_1f_2]} \left(\frac{f_1f_2(u_1k_2,u_2k_1)^2}{(u_1k_2f_1,u_2k_1f_2)^2}\right)^{i\beta}.\end{equation}
	
	In Lemma \ref{varapprox}, we show using Rankin's trick that we may remove the restriction on the number of prime factors of coefficients of the mollifier $M_1...M_l$ in calculating $\lim_{\beta\to 0^+}P_4(\beta)$, with error bounded by $O(\frac{\log q}{\log q_l}\E_{\oplus}[\left|F\right|^2]).$ Note that this is of the same magnitude as the bound in Theorem  \ref{realtwistthm} for the whole twisted moment $\E_{\oplus}\left[\left|LM_1...M_lF\right|^2\right].$
	We want to show $P_4(\beta)$ may be approximated by $\varUpsilon(\beta)$, where 
	
	\begin{equation}\label{wholefsum}
		\varUpsilon(\beta)=\frac{1}{2\beta}\mathlarger{\sum}_{\substack{ (u_1,k_1)=(u_2,k_2)=1}} \left|\tilde{C}_{u_1,k_1}\tilde{C}_{u_2,k_2}\right|\left|\tilde{H}(\beta,u_1,k_1,u_2,k_2)-\tilde{H}(-\beta,u_1,k_1,u_2,k_2)\right|,
	\end{equation} and for  integers $u_1,k_1,u_2,k_2$ with prime factors all at most $q_l$ such that $(u_1,k_1)=(u_2,k_2)=1$, $\tilde{H}(\beta,u_1,k_1,u_2,k_2)$ being defined as:
	\begin{equation}\label{Hunres}			
		\tilde{H}(\beta,u_1,k_1,u_2,k_2)= \mathlarger{\sum}_{\substack{p|f_1f_2\implies p\in (q_{m-1},q_m]}}\frac{ \mu(f_1)\mu(f_2)}{[u_1k_2f_1,u_2k_1f_2]} \left(\frac{f_1f_2(u_1k_2,u_2k_1)^2}{(u_1k_2f_1,u_2k_1f_2)^2}\right)^{i\beta}.\end{equation}
	
	In Lemma \ref{var}, we bound the sum without the restriction on the number of prime factors in the mollifier, and show this may also be bounded as $O\left(\frac{\log q}{\log q_l}\E_{\oplus}[\left|F\right|^2]\right).$
	This will complete the proof of Lemma \ref{S4lem}.
	
	\begin{lem}\label{varapprox}
		Let $P_4(\beta)$ be defined as in \eqref{P4def}, and $\varUpsilon(\beta)$ be as defined in \eqref{wholefsum}, for $F$ an  $l$-sufficient function. Then 
		\begin{equation}	
			P_4(\beta)= \varUpsilon(\beta)+O\left(\frac{\log q}{\log q_l}\E_{\oplus}[\left|F\right|^2]\right)+o_{\beta\to 0^+,F}(1).
		\end{equation}

	\end{lem}
	In order to prove Lemma \ref{varapprox}, it is first necessary to understand the terms defining $\varUpsilon(\beta)$ in Equation \eqref{wholefsum} better, so that we may split it into the contribution from primes in different intervals $(q_{j-1},q_j]$, for $1\le j\le l$. This follows the  proof of the estimate in Lemma \ref{S3approx}, which allowed the restriction on the coefficients of $\tilde{M}$ to be lifted in bounding $P_3.$
	Given positive integers $c_1,c_2,f_1$ and $f_2$ with all their prime factors at most $q_l,$ let \begin{equation}\label{defU}U(c_1,c_2,f_1,f_2)=\frac{f_1f_2(c_1,c_2)^2}{(c_1f_1,c_2f_2)^2}.\end{equation}
	Further define  \begin{equation}\label{theta}
		\Theta_\beta(c_1,c_2)=\mathlarger{\sum}_{p|f_1,f_2\implies p\le q_l}\frac{\mu(f_1)\mu(f_2))}{[c_1f_1, c_2f_2]} \frac{U(c_1,c_2,f_1,f_2)^{i\beta}-U(c_1,c_2,f_1,f_2)^{-i\beta}}{2i\beta}.
	\end{equation}
	Then we see that Equation \eqref{wholefsum} may be rewritten as
	\begin{equation}\label{splitr}	\varUpsilon(\beta)=\mathlarger{\sum}_{u_1,k_1,u_2,k_2} \left|\tilde{C}_{u_1,k_1}\tilde{C}_{u_2,k_2}\right|\left|\Theta_\beta(u_1k_2,u_2k_1)\right|.
	\end{equation}
	We will show that, most tuples $(u_1,k_1,u_2,k_2)$ don't contribute to the summand in \eqref{splitr}.
	Indeed, we will show only the diagonal terms $u_1k_2=u_2k_1$, and those tuples  where the prime factorisations of $u_1k_2$ and $u_2k_1$ differ at a single prime, have non-zero contribution to \eqref{splitr} as $\beta\to0.$
	We will bound both of these contributions, and show that, summing over $j$, both sets of contributions can be bounded by the  desired bound
	$\ll  \E_{\oplus}\left[\left|F\right|^2\right] \times\frac{\log q}{\log q_l}$.
	\begin{lem}\label{P4terms}
		Let $c_1$ and $c_2$ be integers with all their prime factors at most $q_l$.
		\begin{itemize}
			\item If $c_1=c_2$, then \begin{equation}\label{samelemm}
				\Theta_\beta(c_1,c_1)= \frac{1}{c_1}\prod_{p\le q_l}\left(1-\frac{1}{p}\right) \mathlarger{\sum}_{p'\le q_l} \left(\frac{2\log p'}{p'-1}\right)+o_{\beta\to 0^+,F}(1).
			\end{equation}
			\item If there exists a unique prime $p^\ast\le q_l$ such that $v_p(c_1)\ne v_p(c_2)\iff p=p^\ast,$ then \begin{equation}\label{diflem}
				\Theta_\beta(c_1,c_2)=\frac{(p^\ast-1)\log p^\ast}{(p^\ast)^{\max\{v_{p^\ast}(c_1),v_{p^\ast}(c_2)\}+1}} \prod_{\substack{p\le q_l\\ p\ne p^\ast}} \frac{p-1}{p^{\max\{v_p(c_1), v_p(c_2)\}+1}}+o_{\beta\to 0^+,F}(1).	
			\end{equation}
			\item If the prime factorisations of $c_1$ and $c_2$ differ in more than one prime, then \begin{equation}\label{0lem}
				\lim_{\beta\to 0^+} \Theta_\beta(c_1,c_2)=0.		
			\end{equation}
		\end{itemize}
	\end{lem}

	\begin{proof}
		
		Define
		\begin{equation}\label{dephi}
			\Phi_{\beta}(c_1,c_2)=\mathlarger{\sum}_{p|f_1,f_2\implies p\le q_l}\frac{\mu(f_1)\mu(f_2))}{[c_1f_1, c_2f_2]} U(c_1,c_2,f_1,f_2)^{i\beta}.
		\end{equation}
		Then 
		\begin{equation}\label{thetabeta}
			\Theta_\beta(c_1,c_2)=\frac{\Phi_{\beta}-\Phi_{-\beta}}{2i\beta},
		\end{equation}
		and we can calculate $\Phi_{\beta}(c_1,c_2)$ as the sum of functions doubly multiplicative in  $f_1$ and $f_2.$
		
		We first calculate the diagonal terms $\Theta(c_1,c_1)$.
		We have \begin{equation}
			\Phi_{\beta}(c_1,c_1)=\frac{1}{c_1}\mathlarger{\sum}_{p|f_1,f_2\implies p\le q_l} \frac{\mu(f_1)\mu(f_2)}{[f_1,f_2]}\left(\frac{f_1f_2}{(f_1,f_2)^2}\right)^{i\beta}.
		\end{equation}
		
		On the diagonal terms, $\Phi_{\beta}(c_1,c_1)$ factors as an Euler product.
		If we write $d_{p,i}=v_p(f_i)$, then we obtain
		\begin{equation}\label{diageuler}
			\Phi_{\beta}(c_1,c_1)=\frac{1}{c_1}\prod_{p\le q_l} \mathlarger{\sum}_{d_{p,1},d_{p,2}\in \{0,1\}} \frac{(-1)^{d_{p,1}+d_{p,2}}}{p^{\max\{d_{p,1},d_{p,2}\}}}(p^{k_1+k_2-2\min\{k_1,k_2\}})^{i\beta}.\end{equation}
		We can write this as a product over primes in the interval,
		\begin{equation}\label{AA}
			\prod_{p\le q_l} A_{p,\beta}, \qquad  A_{p,\beta}=1+\frac{1}{p}-\frac{2p^{i\beta}}{p}.
		\end{equation}
		Hence, \begin{equation}\label{nosum}
			\Theta_\beta(c_1,c_1)=\frac{1}{2i\beta}(	\Phi_{\beta}(c_1,c_1)-	\Phi_{-\beta}(c_1,c_1))=\frac{\Phi_{-\beta}(c_1,c_1)}{2i\beta}\left(\frac{\Phi_{\beta}(c_1,c_1)}{\Phi_{-\beta}(c_1,c_1)}-1\right).
		\end{equation}
		Each term in the Euler product for $\frac{\Phi_{\beta}}{\Phi_{-\beta}}$ is $1+O(\beta)$, hence we may express \eqref{nosum} as
		\begin{equation}\label{betato0}
			\Theta_\beta(c_1,c_1)=(\Phi_{0}(c_1,c_1)+O(\beta))\mathlarger{\sum}_{p\le q_l} \frac{1}{2i\beta}\left(\frac{A_{p,\beta}}{A_{p,-\beta}}-1\right).
		\end{equation}
		But setting $\beta=0$ in \eqref{AA} shows that:
		\begin{equation}
			\Phi_{0}(c_1,c_1)=\frac{1}{c_1}\prod_{p\le q_l}\left(1-\frac{1}{p}\right).
		\end{equation}
		For each prime $p\le q_l,$
		\begin{equation}
			A_{p,\beta}=\left(1-\frac{1}{p}\right)\left(1-\frac{2i\beta \log p}{p-1}\right)+O(\beta^2),
		\end{equation}
		hence \begin{equation}
			\frac{A_{p,\beta}}{A_{p,-\beta}}=1-\frac{4i\beta \log p}{p-1}+O(\beta^2).
		\end{equation}
		Hence as $\beta\to 0^+$ in \eqref{betato0}, and we take the sum over
		the contribution of all the primes in $(q_{j-1},q_j]$ to the sum on the right-hand side, we obtain:
		\begin{equation}
			\Theta_\beta(c_1,c_1)= \frac{1}{c_1}\prod_{p\le q_l}\left(1-\frac{1}{p}\right) \mathlarger{\sum}_{p'\le q_l} \left(\frac{2\log p'}{p'-1}\right)+o_{\beta\to 0^+,F}(1).
		\end{equation} 
		This proves Equation \eqref{samelemm}.
		
		It remains to calculate the off-diagonal terms for $\Theta_\beta(c_1,c_2).$
		The Euler product, instead of the diagonal term \eqref{diageuler}, becomes the more complicated formula:
		\begin{equation}\label{defphij}
			\Phi_{\beta}(c_1,c_2)=\prod_{p\in (q_{j-1},q_j]} \xi_\beta(p,c_1,c_2), 
		\end{equation}
		where \begin{equation}
			\xi_\beta(p,c_1,c_2)=
			\mathlarger{\sum}_{d_{p,1},d_{p,2}\in \{0,1\}} \frac{(-1)^{d_{p,1}+d_{p,2}}}{p^{\max\{ v_p(u_1k_2)+d_{p,1}, v_p(u_2k_1)+d_{p,2}\}}}\left(\frac{p^{d_{p,1}+d_{p,2}+2\min\{v_p(u_1k_2), v_p(u_2k_1)\}}}{p^{2\min\{d_{p,1}+v_p(u_1k_2), d_{p,2}+v_p(u_2k_1)\}}}\right)^{i\beta}.
		\end{equation}
		
		This is different to most of the Euler products, since some of the terms may be $O(\beta)$, rather than $1+O(\beta)$.
		We remark that this doesn't occur when splitting $\Phi_\beta(c_1,c_1)$ into the $A_{p,\beta}$ in Equation \eqref{betato0}, because the leading term for $A_{p,\beta}$ is $
		A_{p,0}=1-\frac{1}{p}$ which is non-zero, and hence $
		\frac{A_{p,\beta}}{A_{p,-\beta}}=1+O(\beta)$.
		
		We now evaluate $\xi_\beta(p,c_1,c_2);$ the limiting value will depend on whether or not $v_p(c_1)=v_p(c_2).$
		If $v_p(c_1)=v_p(c_2)$, then \begin{equation}
			\xi_\beta(p,c_1,c_2)=\frac{1-2p^{i\beta-1}+p^{-1}}{p^{\max\{v_p(c_1),v_p(c_2k)\}}}.\end{equation}
		
		We expand this to first order terms in $\beta$ to show:
		\begin{equation}\label{samterm}
			\xi_\beta(p,c_1,c_2)=\frac{p-1}{p^{\max\{v_p(c_1), v_p(c_2)\}+1}}+O(\beta).
		\end{equation}	
		If $v_p(c_1)\ne v_p(c_2)$, then
		\begin{equation}
			\xi_\beta(p,c_1,c_2)=\frac{1-p^{i\beta-1}-p^{-i\beta}+p^{-1}}{p^{\max\{v_p(c_1),v_p(c_2)\}}}.
		\end{equation}
		We expand to terms of order $\beta^2$ to show:
		\begin{equation}\label{difterm}
			\xi_\beta(p,c_1,c_2)=\frac{i(p-1)\beta \log p}{p^{\max\{v_p(c_1), v_p(c_2)\}+1}}\times (1+O(\beta)).
		\end{equation}
		
		From \eqref{difterm}, we see that if $v_p(c_1)\ne v_p(c_2)$ for more than one prime $p\le q_l$,  then substituting into \eqref{defphij}, we will have $
		\Phi_{\pm\beta}=O(\beta^2)$,
		since more than one term in the product will be $O(\beta)$.
		Substituting the values for $\Phi_{\pm\beta}$ into  \eqref{thetabeta}, we see that $
		\Theta_\beta(c_1,c_2)=O(\beta)$,
		meaning that the coefficient vanishes in the limit as $\beta$ approaches zero. This proves Equation \eqref{0lem}.
		
		It remains to consider the case where the prime factorisations differ in exactly one prime.
		If we have exactly one prime $p^\ast\le q_l$ such that $v_{p^\ast}(c_1)\ne v_{p^\ast}(c_2)$,
		then  \eqref{difterm} gives expressions for $\xi_\beta(p^\ast,c_1,c_2)$and $\xi_{-\beta}(p^\ast,c_1,c_2)$, while \eqref{samterm} gives expressions for $\xi_\beta(p,c_1,c_2)$ and $\xi_{-\beta}(p,c_1,c_2)$ for all other primes $p\le q_l$ . Substituting these expressions into \eqref{thetabeta}, we see:
		\begin{equation}\label{nomult}
			\Theta_\beta(c_1,c_2)= \frac{2i(p^\ast-1)\beta\log p^\ast}{(p^\ast)^{\max\{v_{p^\ast}(c_1),v_{p^\ast}(c_2)\}+1}} \prod_{\substack{p\le q_l\\ p\ne p^\ast}} \frac{p-1}{p^{\max\{v_p(c_1), v_p(c_2)\}+1}}\times (1+O(\beta))\times \frac{1}{2i\beta}.
		\end{equation}
		This completes the proof of Equation \eqref{diflem}.
	\end{proof}
	
	\begin{proof}[Proof of Lemma \ref{varapprox}]
		From Lemma \ref{P4terms}, we have $\lim_{\beta\to 0^+} \Theta_{\beta}(u_1k_2,u_2k_1)$ is a non-negative real number for all choices of $u_1,k_1,u_2$ and $k_1$.
		Hence, we may remove the final modulus signs in Equation \eqref{splitr} 
		to write
		\begin{equation}\label{splitr2}	\varUpsilon(\beta)=\mathlarger{\sum}_{u_1,k_1,u_2,k_2} \left|\tilde{C}_{u_1,k_1}\tilde{C}_{u_2,k_2}\right|\Theta_\beta(u_1k_2,u_2k_1)+o_{\beta\to 0^+,F}(1).\end{equation}
		
		We now split \eqref{splitr2} into the contribution from primes in different intervals. Given $1\le j\le l$, define \begin{equation}
			\psi(j,\beta)=\mathlarger{\sum}_{m_1,r_1,m_2,r_2} \left|C_{m_1,r_1,B^{(j)}, K_{j}}C_{m_2,r_2,B^{(j)}, K_{j}}\right|\Phi_\beta(m_1r_2,m_2r_1) .
		\end{equation} 
		We may decompose a tuple $(u_1,k_1,u_2,k_2)$, such that $
		p|u_1k_1u_2k_2\implies p \in (q_{j-1},q_j]$ into its factors from each interval.
		If \begin{equation}
			u_1,u_2,k_1,k_2=\prod^l_{j=1}m_{j,1}, \prod^l_{j=1}m_{j,2},\prod^l_{j=1}r_{j,1},\prod^l_{j=1}r_{j,2},
		\end{equation}
		where $
		p|m_{j,1}m_{j,2}r_{j,1}r_{j,2}\implies p\in  (q_{j-1},q_j]$,
		then we may write
		\eqref{splitr2} as:
		\begin{equation}\label{splitr3}	\varUpsilon(\beta)=\prod^l_{j=1}\psi(j,\beta)+o_{\beta\to 0^+,F}(1).\end{equation}
		For any $1\le j\le l$, upon  setting  $i=j$ we have $\psi(j,0)$ is the value of the expression in \eqref{unrestmain}, which is evaluated in Lemma \ref{S3eval} to be the non-zero quantity, $\E_{\oplus}\left[\left|F_j\right|^2\right] \prod_{p\in (q_{j-1},q_j]} (1-\frac{1}{p}).$ Hence, using \eqref{splitr3}, we may write
		\begin{align}\label{newUp}
			\varUpsilon(\beta)=\sum^l_{j=1}& \left(\frac{\psi(j,\beta)-\psi(j,-\beta)}{2i\beta}\right)\prod^l_{\lambda=1,\lambda\ne j} \E_{\oplus}\left[\left|F_\lambda\right|^2\right] \prod_{p\in (q_{\lambda-1},q_\lambda]} (1-\frac{1}{p})	\\\nonumber&+o_{\beta\to 0^+,F}(1).
		\end{align} 
		
		We must bound the contribution to \eqref{newUp} of the extra terms in $\tilde{M}$ both with  $\max\{\Omega_j(f_1),\Omega_j(f_2)\}>10(n_j-n_{j-1})^{10^5}$  for each $j$, and the terms  	with $\max\{\Omega_\lambda(f_1),\Omega_\lambda(f_2)\}>10(n_\lambda-n_{\lambda-1})^{10^5}$ for each $\lambda\ne j$.
		In order to bound such terms, we first need to understand the expression related to numbers with prime factors all lying in the interval $(q_{j-1},q_j].$
		\begin{prop}\label{logjterm}
			Let $1\le j\le l$  and $\psi(j,\beta)$ be as above. Then  
			\begin{equation}
				\label{eqn: logjterm}
				\left(\frac{\psi(j,\beta)-\psi(j,-\beta)}{2i\beta \psi(j,\beta)}\right)\ll (\log q_j-\log q_{j-1})\E_{\oplus}\left[\left|F_j\right|^2\right] \prod_{p\in (q_{j-1},q_j]} (1-\frac{1}{p})+o_{\beta\to0^+,F}(1).\end{equation}
		\end{prop}
		\begin{proof}
			We define a restriction of the function $\Theta_j$ defined in \eqref{theta} to numbers with all their prime factors lying in the interval $(q_{j-1},q_j].$
			Given positive integers $c_1$ and $c_2$ with all their prime factors lying in the interval $(q_{j-1},q_j]$,
			define  \begin{equation}\label{thetalocal}
				\Theta_{j,\beta}(c_1,c_2)=\mathlarger{\sum}_{p|f_1,f_2\implies p\in (q_{j-1},q_j]}\frac{\mu(f_1)\mu(f_2))}{[c_1f_1, c_2f_2]} \frac{U(c_1,c_2,f_1,f_2)^{i\beta}-U(c_1,c_2,f_1,f_2)^{-i\beta}}{2i\beta},
			\end{equation}
			where $U$ is the function defined in \eqref{defU}.
			Then we see that the left-hand side of \eqref{eqn: logjterm} is
			\begin{equation}\label{splitj}	\mathlarger{\sum}_{m_{j,1},r_{j,1},m_{j,2},r_{j,2}} \left|C_{m_{j,1},r_{j,1}, B^{(j)}, K_{j}}C_{m_{j,2},r_{j,2}, B^{(j)}, K_{j}}\right|\left|\Theta_{j,\beta}(m_{j,1}r_{j,2},m_{j,2}r_{j,2})\right|.
			\end{equation}
			Following the proof of Lemma \ref{P4terms}, restricted to primes in the interval $(q_{j-1},q_j]$, we see that, most tuples $(m_{j,1},r_{j,1}m_{j,2},r_{j,2})$ don't contribute to the summand in \eqref{splitr}.
			We have the same result with the condition $p\leq q_l$ replaced by $p\in (q_{j-1},q_j]$.
			
			We use this in Equation \eqref{splitj}, and notice that on the support of the coefficients $C_{m_{j,1},r_{j,1}, B^{(j)}, K_{j}}C_{m_{j,2},r_{j,2}, B^{(j)}, K_{j}}$ we have $m_{j,1}r_{j,2}=m_{j,2}r_{j,1}\iff m_{j,1}=m_{j,2},r_{j,1}=r_{j,2}$. Hence the above equals
			\begin{align}\label{splitjbe}
				&\prod_{p\in(q_{j-1},q_j]}\left(1-\frac{1}{p}\right)\left(\mathlarger{\sum}_{m_{j,1},r_{j,1}} \frac{\left|C_{m_{j,1},r_{j,1}, B^{(j)}, K_{j}}\right|^2}{m_{j,1}r_{j,1}}\hspace{0.5cm}+\right.\\\nonumber&\left.\mathlarger{\sum}\limits_{p^{\ast}\in (q_{j-1}, q_j]}\log(p^\ast)\mathlarger{\sum}_{v_{p}(m_{j,1}r_{j,2})\ne v_{p}(m_{j,2}r_{j,1})\iff p=p^\ast} \frac{\left|C_{m_{j,1},r_{j,1}, B^{(j)}, K_{j}}C_{m_{j,2},r_{j,2}, B^{(j)}, K_{j}}\right|}{	\left[m_{j,1}r_{j,2}, m_{j,2}r_{j,1}\right]}
				\right).
			\end{align}
			The diagonal term corresponding to the first sum in $\eqref{splitj}$ may be handled directly using Lemma \ref{Realexp}, to show it contributes:
			$ \prod_{p\in(q_{j-1},q_j]}\left(1-\frac{1}{p}\right) \E_{\oplus}\left[\left|F_j\right|^2\right]$.
			It remains to show the off-diagonal terms from the second sum can be bounded by $$(\log q_j-\log q_{j-1}) \prod_{p\in(q_{j-1},q_j]}\left(1-\frac{1}{p}\right) \E_{\oplus}\left[\left|F_j\right|^2\right].$$
			
			The contribution of the second sum to \eqref{splitj} may be bounded as:
			\begin{equation}\label{onep}
				\prod_{p\in(q_{j-1},q_j]}\left(1-\frac{1}{p}\right)\left(\mathlarger{\sum}_{p^{\ast}\in (q_{j-1}, q_j]}\log p^\ast\mathlarger{\sum}_{v_{p}(m_{j,1}r_{j,2})\ne v_{p}(m_{j,2}r_{j,1})\iff p=p^\ast} \frac{\left|C_{m_{j,1},r_{j,1}, B^{(j)}, K_{j}}C_{m_{j,2},r_{j,2}, B^{(j)}, K_{j}}\right|}{	\left[m_{j,1}r_{j,2}, m_{j,2}r_{j,1}\right]}
				\right).
			\end{equation}
			Applying Cauchy-Schwarz and using the support of the coefficients $C_{m_{j},r_{j}, B^{(j)}, K_{j}}$ shows that \eqref{onep} may be bounded as:
			\begin{align}\label{onepcs}
				\prod_{p\in(q_{j-1},q_j]}&\left(1-\frac{1}{p}\right)\left(\mathlarger{\sum}_{\substack{p|m_{j,1}r_{j,1}\implies p \in (q_{j-1},q_j]\\ (m_{j,1},r_{j,1})=1}} \left|C_{m_{j,1},r_{j,1}, B^{(j)}, K_{j}}\right|^2\right.\\\nonumber&\left.\mathlarger{\sum}_{\substack{p^{\ast}\in(q_{j-1},q_j]\\(m_{j,2},r_{j,2})=1\\v_{p}(m_{j,1}r_{j,2})\ne v_{p}(m_{j,2}r_{j,1})\iff p=p^\ast}}\log p^\ast\frac{1}	{\left[m_{j,1}r_{j,2}, m_{j,2}r_{j,1}\right]}
				\right).
			\end{align}
			Considering a fixed prime $p^\ast$, the values of coprime pairs of integers  $m_{j,1},r_{j,1}$ and  $m_{j,2},r_{j,2}$ with $v_{p}(m_{j,1}r_{j,2})\ne v_{p}(m_{j,2}r_{j,1})\iff p=p^\ast$  may be expressed as:
			\begin{equation}m_{j,1}=m' (p^\ast)^{d_1},\quad r_{j,1}=r' (p^\ast)^{d_2}, \quad m_{j,2}=m' (p^\ast)^{d_3}, \quad r_{j,2}=r' (p^\ast)^{d_4},
			\end{equation}
			where $p^\ast\nmid m'r'$ and $d_1,d_2,d_3,d_4$ are non-negative integers. The conditions:\\ $(m_{j,1},r_{j,1})=(m_{j,2},r_{j,2})=1$ and $v_{p}(m_{j,1}r_{j,2})\ne v_{p}(m_{j,2}r_{j,1})\iff p=p^\ast$
			mean \begin{equation} \label{minn}
				\min\{d_1,d_2\}=\min\{d_3,d_4\}=0
			\end{equation}
			and 
			\begin{equation}\label{nrel}
				d_1+d_4=d_3+d_2+n,	\end{equation}  
			where $n$ is a 	non-zero integer.
			Hence, using the restriction  the inner sum over $(m_{j,2},r_{j,2})$ in \eqref{onepcs} may be expressed as \begin{equation}\label{innern}
				\mathlarger{\sum}\limits_{n\ne 0}\frac{\log p^\ast}{m_{j,1}r_{j,1}}(p^\ast)^{d_1+d_2-\max\{d_1+d_4,d_2+d_3\}}.\end{equation}
			We handle the terms from positive $n$  and negative $n$  separately. We first consider positive $n$; the case of negative $n$ follows by symmetry.
			
			For positive $n$, we have  \begin{equation}\label{posn}
				d_1+d_2-\max\{d_1+d_4,d_2+d_3\}=d_2-d_4-n,
			\end{equation} so that the contribution to the inner sum in Equation \eqref{onepcs} is a geometric series.
			In the case $d_2>0$ we have $d_1=0$ by Equation \eqref{minn}, so that Equation \eqref{nrel} yields that \begin{equation}\label{drel}d_4\ge d_2.\end{equation}
			There are relatively few choices of $p^\ast$ for which $d_2\ne 0$ since $d_2$ is constrained to have few prime factors.
			Note that  for fixed $d_1,d_2$ and $d_4$, $d_3$ is determined uniquely by Equation \eqref{nrel}. The two cases to consider are when $d_2=0$ and $d_2>0$.
			In the case $d_2=0$,  we see the contribution to \eqref{innern} of positive $n$ is hence
			\begin{equation}\label{d0}
				\ll \mathlarger{\sum}\limits_{p^\ast\in (q_{j-1},q_j]} \log p^\ast\mathlarger{\sum}\limits_{n>0}\frac{1}{m_{j,1}r_{j,1}(p^\ast)^n}
				\ll \frac{1}{m_{j,1}r_{j,1}}\mathlarger{\sum}\limits_{p^\ast\in (q_{j-1},q_j]}\frac{\log p^\ast}{p^\ast}.
			\end{equation}
			Using the Prime Number Theorem, the bound above is
			\begin{equation}\label{d=0}
				\ll  \frac{\log q_j-\log q_{j-1}}{m_{j,1}r_{j,1}}.
			\end{equation}
			
			For the case $d_2>0$, there are at most $10(n_j-n_{j-1})^{10^4}$ choices of $(p^\ast)$ for which $d_2>0$.
			By \eqref{drel}, they each have contribution
			\begin{equation}
				\ll \mathlarger{\sum}\limits_{n>0}\frac{\log p^\ast}{m_{j,1}r_{j,1}(p^\ast)^n}
				\ll \frac{\log q_{j-1}}{m_{j,1}r_{j,1}q_{j-1}} ,
			\end{equation}
			so that the total contribution to \eqref{innern} is \begin{equation}\label{d>0}
				\ll \frac{(n_j-n_{j-1})^{10^4}\log q_{j-1}}{m_{j,1}r_{j,1}q_{j-1}}\ll \frac{1}{m_{j,1}r_{j,1}}.
			\end{equation} 
			Combining the bounds in Equations \eqref{d=0} and \eqref{d>0} we see the total contribution to \eqref{innern} for positive $n$ is 
			\begin{equation}\label{dg}
				\ll\frac{\log q_j-\log q_{j-1}}{m_{j,1}r_{j,1}}.
			\end{equation}
			
			By symmetry, the contribution to \eqref{innern} for negative $n$ is also
			\begin{equation}
				\ll\frac{\log q_j-\log q_{j-1}}{m_{j,1}r_{j,1}}.
			\end{equation}
			Substituting this bound into \eqref{onepcs}, we have that the left-hand side of \eqref{eqn: logjterm} is
			\begin{equation}
				\ll (\log q_j-\log q_{j-1})\prod_{p\in(q_{j-1},q_j]}\left(1-\frac{1}{p}\right)\left(\mathlarger{\sum}_{\substack{ (m_{j,1},r_{j,1})=1}} \frac{\left|C_{m_{j,1},r_{j,1}, B^{(j)}, K_{j}}\right|^2}{m_{j,1}r_{j,1}}
				\right).
			\end{equation}
			Using Lemma \ref{Realexp} to handle the sum shows that the right-hand side is
			\begin{equation}
				(\log q_j-\log q_{j-1})\prod_{p\in(q_{j-1},q_j]}\left(1-\frac{1}{p}\right)\E_{\oplus}\left[F_j^2\right].
			\end{equation}
			This completes the proof of Proposition \ref{logjterm}.
		\end{proof}
		We are now in the position to show that we may lift the restriction on the prime factors of coefficients of $\tilde{M}$ in calculating $\lim_{\beta\to0^+}P_4(\beta).$	
		By Rankin's trick, the contribution to  the term $\left(\frac{\psi(j,\beta)-\psi(j,-\beta)}{2i\beta}\right)$ in the summand in \eqref{newUp} of those $f_1$ with $\Omega_j(f_1)>10(n_j-n_{j-1})^{10^5}$ for some $1\le j\le l$ may be bounded for any $\rho >0$ by:	
		\begin{align}\label{CherH}
			e^{-10\rho(n_j-n_{j-1})^{10^5}}&\mathlarger{\sum}_{m_1,r_1,m_2,r_2
			} \left|C_{m_1,r_1,B^{(j)}, K_{j}}C_{m_2,r_2,B^{(j)}, K_{j}}\right|\cdot \\\nonumber&\\\nonumber&\mathlarger{\sum}_{p|f_1f_2\implies p\in (q_{j-1},q_j]}e^{\rho\Omega_j(f_1)}\frac{ \mu^2(f_1)\mu^2(f_2)}{[m_1r_2f_1,m_2r_1f_2]} \left|\log\left(\frac{f_1f_2(m_1r_2,m_2r_1)^2}{(m_1r_2f_1,m_2r_1f_2)^2}\right)\right|+o_{\beta\to 0^+,F}(1).
		\end{align}
		The restriction on the number of prime factors in the interval $(q_{j-1},q_j]$ in \eqref{CherH} is inherited from the support of $C_{m,r,B^{(j)}, K_{j}}.$	
		But since \begin{equation}
			1\le\frac{ f_1f_2(m_1r_2,m_2r_1)^2}{(m_1r_2f_1,m_2r_1f_2)^2}\le f_1f_2,
		\end{equation}
		we see:
		\begin{equation}\left|\log\left(\frac{f_1f_2(m_1r_2,m_2r_1)^2}{(m_1r_2f_1,m_2r_1f_2)^2}\right)\right|\le \log (f_1f_2).\end{equation}
		Hence, the right-hand side in \eqref{CherH} may be bounded by:	
		
		\begin{align}\label{CherH2}
			e^{-10\rho(n_j-n_{j-1})^{10^5}}&\mathlarger{\sum}_{m_1,r_1,m_2,r_2
			} \left|C_{m_1,r_1,B^{(j)}, K_{j}}C_{m_2,r_2,B^{(j)}, K_{j}}\right|\cdot \\\nonumber&\\\nonumber&\mathlarger{\sum}_{p|f_1f_2\implies p\in (q_{j-1},q_j]}e^{\rho\Omega_j(f_1)}\frac{ \mu^2(f_1)\mu^2(f_2)}{[m_1r_2f_1,m_2r_1f_2]} \log\left(f_1f_2\right)+o_{\beta\to 0^+,F}(1) .
		\end{align}
		For each tuple $(m_1,r_1,m_2,r_2)$, define
		\begin{equation} \label{dj}
			D_j(m_1,r_1,m_2,r_2)=\mathlarger{\sum}_{p|f_1f_2\implies p\in (q_{j-1},q_j]}e^{\rho\Omega_j(f_1)}
			\frac{\mu^2(f_1)\mu^2({f_2)}}{[m_1r_2f_1,m_2r_1f_2]}\log\left(f_1f_2\right)
		\end{equation}
		to be the coefficient multiplying $\left| C_{u_1,k_1,K_{(j)}, B^{(j)}}C_{u_2,k_2,K_{(j)}, B^{(j)}}\right|$ in  \eqref{CherH2}, arising from the inner sum over values of $f_1$ and $f_2$. 
		
		The term $\log(f_1f_2)$ in \eqref{dj} may be written by factorising $f_1f_2$ into its prime factors, so that the sum  can be rewritten as a sum over logarithms of primes in the interval. Indeed, we have \begin{equation}
			\log(f_1f_2)=\mathlarger{\sum}_{p|f_1f_2} (\log p) v_p(f_1f_2).
		\end{equation}
		If we write $d_{p,i}=v_p(f_i)$,
		and observe that the exponents can be at most $1$ to contribute to the sum, we see:
		
		\begin{equation}\label{Dnew}
			D_j(m_1,r_1,m_2,r_2)\ll \mathlarger{\sum}_{p\in (q_{j-1},q_j]} \log p \mathlarger{\sum}_{d_{p,1},d_{p,2}\in \{0,1\}} e^{\rho \Omega_j(p^{d_{p,1}})}\times \frac{1}{p^{\max \{d_{p,1}+v_p(m_1r_2),d_{p,2}+v_p(m_2r_1)\}}}\times  \prod_{\substack{t\le q_l\\ t\ne p}} h_{t},
		\end{equation}
		where for each prime $t\in (q_{j-1},q_j]$,  we define:	
		\begin{equation}\label{ht}
			h_{t}=\mathlarger{\sum}\limits_{d_{t,1},d_{t,2}\in \{0,1\}}\frac{e^{\rho \Omega_j(t^{d_{t,1}})}}{t^{\max \{d_{t,1}+v_t(m_1r_2), d_{t,2}+v_t(m_2r_1)\}}}.
		\end{equation} 
		The sum  over $d_{p,1}$ and $d_{p,2}$ is at most: \begin{equation}\label{pb}
			\frac{1}{p^{\max \{v_p(m_1r_2),v_p(m_2r_1)\}}}\times \begin{cases} \frac{100 e^{\rho}\log p}{p} & \text{if } v_p(u_1k_2)=v_p(u_2k_1),\\
				100  e^{\rho}\log p & \text{otherwise}\end{cases}.
		\end{equation}
		
		The sum over $d_{t,1}$ and $d_{t,2}$ for $t$ a prime in the interval not equal to $p$ in \eqref{ht}, which we have denoted by $h_{t}$, is at most: \begin{equation} \label{tb}
			\frac{1}{t^{\max \{v_t(m_1r_2),v_t(m_2r_1)\}}}\times \begin{cases}1+ \frac{100 e^{\rho}}{t} & \text{if } v_t(m_1r_2)=v_t(r_2m_1)\\
				100  e^{\rho} & \text{otherwise}\end{cases}.
		\end{equation}
		
		Hence, substituting the bounds from \eqref{pb} and \eqref{tb} into \eqref{Dnew}, we may bound $	D_{j}(m_1,r_1,m_2,r_2)$ by:
		\begin{align}
			D_{j}(m_1,r_1,m_2,r_2)	\le \frac{1}{[m_1r_2,m_2r_1]}&\times  \mathlarger{\sum}_p \begin{cases} \frac{100 e^{\rho}\log p}{p} & \text{if } v_p(m_1r_2)=v_p(m_2r_1)\\
				100  e^{\rho}\log p & \text{otherwise}\end{cases}\\\nonumber &\times \prod_{t\ne p} \begin{cases}1+ \frac{100 e^{\rho}}{t} & \text{if } v_t(m_1r_2)=v_t(m_2r_1)\\
				100  e^{\rho} & \text{otherwise}\end{cases}.
		\end{align}
		We remove the restriction on $t\ne p$ to obtain that:
		\begin{align}	\label{Dunres}
			D_{j}(m_1,r_1,m_2,r_2)\le \frac{1}{[m_1r_2,m_2r_1]}&\times  \mathlarger{\sum}_p \begin{cases} \frac{100 e^{\rho\Omega_j(p)}\log p}{p} & \text{if } v_p(m_1r_2)=v_p(m_2r_1)\\
				100  e^{\rho\Omega_j(p)}\log p & \text{otherwise}\end{cases}\\\nonumber &\times \prod_{t} \begin{cases}1+ \frac{100 e^{\rho\Omega_j(t)}}{t} & \text{if } v_t(m_1r_2)=v_t(m_2r_1)\\
				100  e^{\rho\Omega_j(t)} & \text{otherwise}\end{cases}.
		\end{align}
		Observe that the above expression factors into a sum over primes $p$ in the interval, and a product over primes $t$. We bound the contribution of both the sum over $p$ and the product over $t$.
		The contribution of those $p$ with $v_p(m_1r_2)=v_p(m_2r_1)$ to the sum in \eqref{Dunres} may be bounded as:
		\begin{equation}
			\le \mathlarger{\sum}_{p\in (q_{j-1}, q_j] }\frac{100 e^\rho \log p}{ p}
			\ll e^\rho e^{(n_j-n_{j-1})},
		\end{equation}
		where we used the Prime Number Theorem to get the final bound. 
		
		Due to the  restriction on the number of prime factors of $m_1,m_2,r_1$ and $r_2$, there are at most $4(n_j-n_{j-1})^{10^4}$ primes $p$ with $v_p(m_1r_2)\ne v_p(m_2r_1)$ in the interval $(q_{j-1},q_j]$.
		Hence, the contribution of those $p$ with $v_p(u_1k_2)\ne v_p(u_2k_1)$, of which to the sum in \eqref{Dunres} is
		\begin{equation}\label{pboundunres}\ll  (n_j-n_{j-1})^{10^4}(\log q_{j}) e^\rho.\end{equation}
		
		We now turn to bound the product over $t$ in \eqref{Dunres}. The contribution of those $t$ for which $v_t(m_1r_2)=v_t(m_2r_1)=0$ to the product in \eqref{Dunres} is
		\begin{equation}\ll  \prod_{t\in (q_{j-1},q_j]} 1+\frac{100e^\rho}{t}\ll e^{100e^\rho(n_j-n_{j-1})} .\end{equation}
		Hence, we may bound the product over $t$ in \eqref{Dunres} as 
		\begin{equation}\label{tbound}\ll e^{100e^\rho(n_j-n_{j-1})}\frac{\alpha(m_1r_2,m_2r_1)}{[m_1r_2,m_2r_1]} ,\end{equation} where $\alpha$ is the multiplicative function defined as in \eqref{alpha}.
		Combining the bound in \eqref{pboundunres} for the sum over $p$ with the bound in \eqref{tbound} for the product over $t$ in \eqref{Dunres}, we obtain:

		\begin{equation} D_j(m_1,r_1,m_2,r_2)\ll  (n_j-n_{j-1})^{10^4}(\log q_j) e^{100 e^\rho(n_j-n_{j-1})} \frac{\alpha(m_1r_2,m_2r_1)}{[m_1r_2,m_2r_1]}.\end{equation} 
		Substituting this bound for the internal sum over $f_1$ into \eqref{CherH}, we see the contribution of the extra values of $f_1$ and $f_2$ to the term $\left(\frac{\psi(j,\beta)-\psi(j,-\beta)}{2i\beta}\right)$ is at most:
		\begin{equation}
			\begin{aligned}\label{concCher}
				&e^{-10\rho(n_j-n_{j-1})^{10^5}}(n_j-n_{j-1})^{10^4}(\log q_j) e^{100 e^\rho(n_j-n_{j-1})} \times \\
				&\hspace{3cm} \mathlarger{\sum}_{m_1,r_1,m_2,r_2}\left| C_{m_1,r_1,K_{(j)}, B^{(j)}}C_{m_2,r_2,K_{(j)}, B^{(j)}}\right|\frac{\alpha(m_1r_2,m_2r_1)}{[m_1r_2,m_2r_1]}+o_{\beta\to 0^+,F}(1).
			\end{aligned} 
		\end{equation}
		proceed as in the proof of Lemma \ref{S3approx}.
		The expression in \eqref{concCher} may be bounded by:
		\begin{equation}
			\begin{aligned}\label{concCher2}
				&e^{-10\rho(n_j-n_{j-1})^{10^5}+100 e^\rho(n_j-n_{j-1})}(n_j-n_{j-1})^{10^4}(\log q_j) \times\\
				&\mathlarger{\sum}_{p|m_2r_2\implies p\in (q_{j-1},q_j]}\left|C_{m_2,r_2,K_{(j)}, B^{(j)}}\right|^2  	\mathlarger{\sum}\limits_{\substack{p|m_1r_1\implies p\in (q_{j-1},q_j]\\ 
						(m_1,r_1)=(m_2,r_2)=1}} \frac{\alpha(m_1r_2,m_2r_1)}{[m_1r_2,m_2r_1]}+o_{\beta\to 0^+,F}(1).\end{aligned}
		\end{equation}
		Using \eqref{innerto2} and \eqref{sumbound} to bound the internal sum over $m_1$ and $r_1$,
		we may bound \eqref{concCher2} as  
		\begin{align}\label{concCher3}
			\ll	e^{-10\rho(n_j-n_{j-1})^{10^5}+400e^\rho(n_i-n_{i-1})+2.10^4 \rho(n_j-n_{j-1})^{10^4+1}}&(n_j-n_{j-1})^{10^4}(\log q_j) e^{100 e^\rho(n_j-n_{j-1})}\\\nonumber\mathlarger{\sum}_{\substack{p|m_2r_2\implies p\in (q_{j-1},q_j]}}\frac{\left| C_{m_2,r_2,K_{(j)}, B^{(j)}}\right|^2}{m_2r_2} +o_{\beta\to 0^+,F}(1).&\end{align}
		Observing that the sum over $m_2,r_2$ is $\E_{\oplus}\left[\left|F_j\right|\right]^2$, we may bound \eqref{concCher3} as
		\begin{equation}\label{concCher4}
			\ll	\log(q)e^{-10\rho(n_j-n_{j-1})^{10^5}+600e^\rho(n_i-n_{i-1})+2.10^4 \rho(n_j-n_{j-1})^{10^4+1}}\E_{\oplus}\left[\left|F_j\right|^2\right] +o_{\beta\to 0^+,F}(1).\end{equation}
		Taking $\rho =1000$, the above is 
		$\ll 	\log(q)e^{-100(n_j-n_{j-1})}\E_{\oplus}\left[\left|F_j\right|^2\right] +o_{\beta\to 0^+,F}(1).$ By symmetry, the contribution of those values of $f_2$ with too many prime factors to the term $\left(\frac{\psi(j,\beta)-\psi(j,-\beta)}{2i\beta}\right)$ is also 	\begin{equation}\ll 	\log(q)e^{-100(n_j-n_{j-1})}\E_{\oplus}\left[\left|F_j\right|^2\right] +o_{\beta\to 0^+,F}(1).\end{equation}  
		Substituting this bound into \eqref{newUp}, the contribution of the the extra terms $f_1$ and $f_2$ to $\left(\frac{\psi(j,\beta)-\psi(j,-\beta)}{2i\beta}\right)$, summed over all $1\le j\le l$,  is \begin{equation}\label{jCher}
			\ll\sum^l_{j=1} \log(q_j)e^{-100(n_j-n_{j-1})}\E_{\oplus}\left[\left|F_j\right|^2\right].\prod^l_{\lambda=1,\lambda\ne j} \E_{\oplus}\left[\left|F_\lambda\right|^2\right] \prod_{p\in (q_{\lambda-1},q_\lambda]} (1-\frac{1}{p})	+o_{\beta\to 0^+,F}(1).
		\end{equation} 
		Using Lemma \ref{Realsplit}, the bound in \eqref{jCher} is $\ll \E_{\oplus}[F^2].$ 
		
		It remains to show the contribution of the values of $f_1$ and $f_2$ with $\max\{\Omega_\lambda(f_1),\Omega_\lambda(f_1)\}>10(n_\lambda-n_{\lambda-1})^{10^5}$ to the expression in \eqref{newUp} for $\lambda\ne j$ can also be bounded by $\E_{\oplus}\left[F^2\right]$.
		By Lemma \ref{S3approx}, we see that the contribution of the extra terms $f_1$ and $f_2$ to $\E\left[F_\lambda^2\right]$ is $O\left(e^{-100(n_\lambda-n_{\lambda-1})}\E_{\oplus}\left[\left|F_\lambda\right|^2\right]\right)$.
		By Proposition \ref{logjterm}, the term associated to the interval $(q_{j-1},q_j]$ in Equation \eqref{newUp} is
		\begin{equation}\left(\frac{\psi(j,\beta)-\psi(j,-\beta)}{2i\beta \psi(j,\beta)}\right)\ll (\log q_j-\log q_{j-1})\E_{\oplus}\left[\left|F_j\right|^2\right] \prod_{p\in (q_{j-1},q_j]} (1-\frac{1}{p})+o_{\beta\to_0^+,F}(1).\end{equation}
		Hence, the contribution of all the extra terms in the intervals with 
		$\lambda\ne j$ to \eqref{newUp} is
		\begin{align}
			\ll &\sum^l_{j=1} (\log q_j-\log q_{j-1})\E_{\oplus}\left[\left|F_j\right|^2\right]\prod_{p\in (q_{j-1},q_j]} (1-\frac{1}{p})\sum^l_{\lambda=1,\lambda\ne j} e^{-100(n_\lambda-n_{\lambda-1})}\E_{\oplus}\left[\left|F_\lambda\right|^2\right] 	\\&\prod_{p\in (q_{\lambda-1},q_\lambda]} (1-\frac{1}{p}) \prod^l_{\mu=1, \mu \ne\lambda,j}\E_{\oplus}\left[\left|F_\mu\right|^2\right] \prod_{p\in (q_{\mu-1},q_\mu]} (1-\frac{1}{p})\nonumber+o_{\beta\to 0^+,F}(1).
		\end{align}
		An application of the Prime Number Theorem and Lemma \ref{Realsplit} shows the above bound is 
		\begin{equation}
			\ll \frac{\E_{\oplus}\left[F^2\right]}{\log q_l}\sum^l_{j=1} (\log q_j-\log q_{j-1})\sum^l_{\lambda=1,\lambda\ne j} e^{-100(n_\lambda-n_{\lambda-1})}+o_{\beta\to 0^+,F}(1).
		\end{equation}
		The sum over $\lambda\ne j$ is bounded uniformly for each $1\le j\le 1$, hence the contribution of the additional primes to the terms $\lambda\ne j$ is $	\ll\E_{\oplus}{\left[F^2\right]}.$ This completes the proof of Lemma \ref{varapprox}.\end{proof}

	We now bound $\varUpsilon(\beta)$.
	\begin{lem}\label{var}
		Let $\varUpsilon(\beta)$ be as defined as in \eqref{wholefsum}. Then 
		\begin{equation}
			\varUpsilon(\beta)\ll  \frac{\log q}{\log q_l}\E_{\oplus}\left[\left|F\right|^2\right] +o_{\beta\to 0^+,F}(1).
		\end{equation}
	\end{lem}
	
	\begin{proof}
		Recall from \eqref{newUp}: 
		\begin{align}\label{newUpcopy}
			\varUpsilon(\beta)=\sum^l_{j=1} \left(\frac{\psi(j,\beta)-\psi(j,-\beta)}{2i\beta}\right)\prod^l_{\lambda=1,\lambda\ne j} \E_{\oplus}\left[\left|F_\lambda\right|^2\right] \prod_{p\in (q_{\lambda-1},q_\lambda]} (1-\frac{1}{p})+o_{\beta\to 0^+,F}(1).
		\end{align}
		By Proposition \ref{logjterm}, we may rewrite \eqref{newUpcopy} as
		\begin{align}\label{newUpcopy2}
			\varUpsilon(\beta)\ll \sum^l_{j=1}& (\log q_j-\log q_{j-1})\E_{\oplus}\left[\left|F_j\right|^2\right] \prod_{p\in (q_{j-1},q_j]} \prod^l_{\lambda=1,\lambda\ne j} \E_{\oplus}\left[\left|F_\lambda\right|^2\right] \prod_{p\in (q_{\lambda-1},q_\lambda]} (1-\frac{1}{p})	\\\nonumber&+o_{\beta\to 0^+,F}(1).
		\end{align}
		Using Lemma \ref{Realsplit} and the Prime Number Theorem yields 
		\begin{equation}\label{newUpcopy3}
			\varUpsilon(\beta)\ll \E_{\oplus}[F^2]\sum^l_{j=1} (\log q_j-\log q_{j-1})\prod_{p\le q_l} (1-\frac{1}{p})	+o_{\beta\to 0^+,F}(1)\ll \E_{\oplus}[F^2]:+o_{\beta\to 0^+, F(1)}.
		\end{equation}
	\end{proof}
	This completes the proof of Lemma \ref{var}. Using Lemma \ref{varapprox} and Lemma \ref{var} we conclude the proof of Proposition \ref{S4lem}, and hence the proof of Theorem \ref{realtwistthm}.
	
	\begin{proof}[Proof of Theorem \ref{twistthm}] 
		This essentially follows from writing $|Q_j|^2=\Re(Q_j^2)+\Re(iQ_j)^2$. Taking $K_j(x)=x$ to be the identity polynomial and for each subset $S\subset\{1,...,l\}$ taking $F_j$ to be $Q_j i^{{\rm 1}(j\in S)}$  in Theorem \ref{realtwistthm}, we obtain  
		\begin{equation}
			\mathlarger{\sum}\limits_{S\subset \{1,...,l\}} \E_{\oplus}\left[\left|LM\right|^2 \prod^l_{j=1}\Re(Q_j i^{{\rm 1}(j\in S)})^2\right]\ll \frac{\log q}{\log q_l}\mathlarger{\sum}\limits_{S\subset \{1,...,l\}} \E_{\oplus}\left[\prod^l_{j=1}\Re(Q_j i^{{\rm 1}(j\in S)})^2\right], 
		\end{equation}
		or  
		\begin{equation}
			\E_{\oplus}\left[\left|LMQ\right|^2\right] \ll \frac{\log q}{\log q_l} \E_{\oplus}\left[\left|Q\right|^2\right] . 
		\end{equation}
	\end{proof}

	\section{Large Deviations}\label{devsec}
	In this section, we prove Theorems \ref{LargeDev} and \ref{largerange}.
	\subsection{Method of proof}The orthogonality relations in the previous section are very similar to the orthogonality relations used in the $t$-aspect in \cite{AB}, and Theorem \ref{twistthm} provides the analogue for the twisted second moment of zeta. Using these points of reference, the method of proof for the large deviations in the $q$-aspect uses a recursive scheme in the same spirit as  \cite{AB}.
	
	We consider the event 
	\begin{equation}
		H=\{\log(L(1/2,\mathfrak{X})>V)\},
	\end{equation}
	where $\mathfrak{X}$ is drawn uniformly at random from the even primitive characters with modulus $q$.
	Since the number of even, primitive characters is asymptotic to $\frac{\varphi(q)}{2}$, we have 
	\begin{equation}\frac{1}{\varphi(q)} \# \{\chi \text{ even,  primitive mod }q: \left|L(1/2,\chi)\right| >V\}\asymp \P(H).\end{equation}
	We partition the event $H$ on {\it good events} pertaining to the partial sums $S_{n_l}$ defined in \eqref{eqn: S} (and \eqref{tildes} for the complex ones). 
	More precisely, we consider the restrictions $\{S_{n_l}\in [L_l, U_l]\}$ where
	\begin{equation}
		L_l=\kappa n_l -\mathbf{s} (\log_{l+2} q) \qquad U_l= \kappa n_l+\mathbf{s} (\log_{l+2} q),
	\end{equation}	
	where $\kappa$ is the {\it slope}
	\begin{equation}
		\kappa=\frac{V}{\log\log q}.
	\end{equation}
	We recall that $V\sim \alpha \log\log q$ for some $0<\alpha<1,$ so that $\kappa\to \alpha.$
	The idea is that if that the ``full sum" $\log\left(L(1/2,\mathfrak{X})\right)$ is close to $\kappa \log\log q$, then the most likely values for $S_{n_l}$ should be $\kappa n_l$. 
	The terms $\pm \mathbf{s} (\log_{l+2} q)$ handles the fluctuations around the linear interpolation.
	To properly mollify  $\log\left(L(1/2,\mathfrak{X})\right)$, more restrictions are needed. 
	We recursively define the following events, where the dependence on $\frac{1}{2}$ and $\mathfrak{X}$ is suppressed for simplicity:\\
	\begin{equation}
		\begin{aligned}
			G_l&=A_l\cap B_l\cap C_l\cap D_l\\
			A_l&= A_{l-1}\cap \left\{\left|\tilde{S}_{n_l}-\tilde{S}_{n_{l-1}}\right| \le \mathscr{A}(n_{l-1}-n_l)\right\}\\
			B_l&=B_{l-1}\cap \left\{S_{n_l}\le U_l\right\}\\
			C_l&=C_{l-1}\cap \left\{S_{n_l}\ge L_l\right\}\\
			D_l&=D_{l-1}\cap \left\{\left|Le^{-S_{n_l}}\right|^2 \le c_l\left|LM_1...M_l\right|^2+e^{-\mathscr{D}(\log\log q-n_{l-1})}\right\}.
		\end{aligned}
	\end{equation}
	where $c_l=\prod^{l}_{j=1}(1+e^{-n_{j-1}})$, and $A_0,B_0,C_0, D_0=\{\text{primitive even characters mod } q\}$ is the full sample space. 
	
	
	
	The following proposition gathers the necessary estimate to bound the probability of $H$ on the good events.

	\begin{prop}
		\label{prop G}
		Let $V= \kappa\log\log q$ with $\frac{1}{\sqrt{\log\log q}}\ll \kappa<1$. 
		Then for some $\delta\gg\kappa$ dependent on $\kappa$, we have
		\begin{enumerate}[(i)]
			\item $\P(H\cap G_1^c )\ll  \frac{e^{-\frac{V^2}{\log\log q}}}{\sqrt{\log\log q}}(\log\log q)^{-\delta}$;
			\item for $1\leq l\leq \mathcal L-1$, $\P(H\cap G_l\cap G_{l+1}^c )\ll  \frac{e^{-\frac{V^2}{\log\log q}}}{\sqrt{\log\log q}}.(\log_{l+2} q)^{-\delta}$;
			\item $\P(H\cap G_{\mathcal L}^c )\ll  \frac{e^{-\frac{V^2}{\log\log q}}}{\sqrt{\log\log q}}$.
		\end{enumerate}
	\end{prop}
	
	%
	\begin{proof}[Proof of Theorem \ref{LargeDev}]
		The proof of Theorem \ref{LargeDev} is completed by writing
		\begin{equation}\label{Hprob}
			\P(H)\le \P(H\cap G_1^c)+\mathlarger{\sum}^{\mathscr{L}-1}_{l=1}\P(H\cap G_l\cap G_{l+1}^c)+P(H\cap G_{\mathscr{L}})
		\end{equation}
		and observing that $\sum^{\mathscr{L}}_{l=1}(\log_{l+2}(q))^{-\delta}=O(1)$ as $q\to \infty$.
	\end{proof}
	
	\begin{proof}[Proof of Theorem \ref{largerange}] We employ the same method to prove Theorem \ref{largerange}, but since we may take $V=o(\log\log q)$, we cannot obtain as sharp bounds as in Theorem \ref{LargeDev}. If $V=o(\log\log q)$ and $N$ is a fixed positive integer, then by the construction of $\mathscr{L}$ we  have $ \log_{\mathscr{L}-n}(q)\gg 1$ as $q\to\infty$  for each $1\le n\le N$. Since $\kappa=o(1),$ we can only obtain $\delta=o(1)$ from these bounds.
		Hence the contribution from $	\P(H\cap G_{\mathscr{L}-n}\cap G_{\mathscr{L}+1-n}^c )$ in Proposition \ref{prop G} may be $O(1).\frac{e^{-\frac{V^2}{\log\log q}}}{\sqrt{\log\log q}}$, and we obtain the weaker bound:
		\begin{equation}
			\P(H)\le \mathscr{L}\times \frac{e^{-\frac{V^2}{\log\log q}}}{\sqrt{\log\log q}}.
		\end{equation}
	\end{proof}
	
	\subsection{Choice of parameters}\label{choiceparameters} In order to ensure the bounds on $\P(H)$ are sharp enough to prove Proposition \ref{prop G},  we require certain bounds on the parameters to ensure certain inequalities below are satisfied. 
	For Equation \eqref{bound3}, we will require:
	\begin{equation}\label{Afix}
		1+\mathbf{s}(\kappa^2-\mathscr{A}^2+2\kappa)<0,
	\end{equation}
	whilst for Equation \eqref{bound4} we will also require $\mathbf{s}$ to satisfy:
	\begin{equation}\label{s2fix}
		\frac{1}{2}+\kappa^2+2(\kappa-1)\mathbf{s}<0.
	\end{equation}
	Equation \eqref{s2fix} forces us to take $\mathbf{s}$ proportional to $\frac{1}{1-\kappa}$. We take
	\begin{equation}
		\label{eqn: s parameter}
		\mathbf{s}=\frac{10^5}{1-\kappa},
	\end{equation} which ensures Equation \eqref{s2fix} holds.
	The factor $10^5$ is to ensure \eqref{L}
	Taking $\mathscr{A}=10^3$ n ensures Equation \eqref{Afix} holds.
	Finally, we take $\mathscr{D}=10^4$
	as in \cite{ABR20}, which gives sufficiently good bounds on $\P(A_{l}\cap D_l^c)$, to apply Lemma \ref{23ABR}.
	
	With these choice of parameters, we see that given $\lambda<1$, there exists some $\epsilon<0$ depending on $\lambda$ such that the difference of left-hand side and the right-hand side of the inequalities in Equations \eqref{Afix} and \eqref{s2fix} are smaller than $\epsilon \kappa$ for all  $0<\kappa\le \lambda$. This enables us to consider the range $\kappa\to 0$ in \ref{largerange}, i.e., $V=o(\log\log q).$ In contrast, the inequalities with the choice of parameters in \cite{AB} do not hold with this uniformity, which prevents them from extending their range of $V$.

	\subsection{Proof of Proposition \ref{prop G}(i)}
	We split the event $H\cap G_1^c$ as
	\begin{equation}
		H\cap G_1^c\subset A_1^c \cup B_1^c \cup (H\cap C_1^c \cup A_1\cap D_1)\cup (D_1^c\cap A_1),\end{equation}
	We bound each of the probabilities of each of the four event on the right-hand side, very similarly as in the proof of Proposition 2.1 in \cite{AB}.
	We give the details for $A_1^c$ and $D_1^c\cap A_1$. The estimates involving $B_1$ and $C_1$ are similar to the ones done in the next section, but simpler. 
	
	Clearly, $\P(A_1^c)=\P\left(|\tilde{S}_{n_1}|>\mathscr{A}n_1\right)$. Thus, using 
	\eqref{softbound} with $k=\lceil 2\mathscr{A}^2n_1\rceil$, we obtain
	\begin{equation}
		\P(A_1^c)	\ll  \sqrt{n_1}\exp(-\mathscr{A}^2 n_1)\ll \frac{e^{-\kappa^2\log\log q}}{\sqrt{\log\log q}}(\log\log q)^{-\delta},
	\end{equation}
	for some $\delta>0$.
	
	%

	
	For $\P(D_1^c\cap A_1)$, we first require a lemma, which allows the restriction $A_1$ to force $L$ to be large unless $D_1$ is satisfied.
	\begin{lem}\label{23ABR}
		Let $l\ge 1$ and suppose  $|\tilde{S}_{n_l}-\tilde{S}_{n_{l-1}}|\le 10^3(n_l-n_{l-1})$.
		Then 
		\begin{equation}\label{mollapprox}
			e^{-(S_{n_l}-S_{n_{l-1}})}\le (1+e^{-n_{l-1}})\left|M_l\right|+e^{-10^5(n_l-n_{l-1})}\end{equation}
	\end{lem}
	\begin{proof}
		This follows from the proof of Lemma 23 in \cite{ABR20} using the identity
		$$
		\prod_{p\in (q_{l-1},q_{l}]}(1-\chi(p)p^{-s})=e^{-(\tilde S_{n_l}-\tilde S_{n_{l-1}}) + R}
		$$
		where $R=\sum_{k\geq 3}\sum_{p\in (q_{l-1},q_{l} ]}\frac{1}{k}p^{-ks}$.
	\end{proof}
	
	The event $A_1\cap \{|L|\le (\log q)^2\}$ is contained in $A_1\cap D_1$. This is because, on $A_1$, one can apply Lemma \ref{23ABR} to get
	$$
	|Le^{-S_{n_1}}|\leq 2|LM_1|+|L|e^{-10^5n_1}\leq 2|LM_1|+(\log q)^2 e^{-10^5n_1}.
	$$
	Hence, $\P(D_1^c\cap A_1)\le \P(L> (\log q)^2)$.
	By Markov's inequality this probability is
	\begin{equation}\label{not D}
		\le (\log q)^{-4} \E_{+}[|L|^2]
		\ll  \frac{e^{-\kappa^2\log q}}{\sqrt{\log\log q}} e^{-\log\log q}.
	\end{equation} 
	Here we used that \begin{equation}
		\label{eqn: L2}
		\E_{+}[|L|^2]\sim \log q,
	\end{equation}
	which is a standard bound and may be seen by taking $M=1$ in \eqref{Qformula}.
	
	\subsection{Proof of Proposition \ref{prop G}(ii)}
	
	As for the proof of Proposition \ref{prop G}(i), we use a union bound on the event $H\cap G_l\cap G_{l+1}^c$ and bound each probability in the union bound individually.
	
	We have the containment: 
	\begin{equation}\label{2.2decomp}H\cap G_l\cap G_{l+1}^c\subset (A_{l+1}^c\cap G_l)\cup (B_{l+1}^c\cup G_l)\cap (H\cap C_{l+1}^c\cap A_{l+1}\cap D_{l+1}\cap G_l)\cup (D_{l+1}^c\cap A_{l+1}\cap G_l),\end{equation}
	In order to bound the probabilities of each event, we need two important lemmas that express the probabilities of the events involving the partial sums in terms the random model where the characters $\chi$ are replaced by the random phases $X$ given in \eqref{eqn: X}. 
	We first prove an analogue of Lemma 2.6 in \cite{AB}, which allows us to bound the twisted moment with the truncated sums bounded above and below by conditioning on $B_l\cap C_l.$
	\begin{lem}\label{2.6}
		Let $1\leq l\leq \mathcal L$.
		Let $Q_l=Q_l(\chi)$ be a Dirichlet polynomial of length $N\le q^{\frac{1}{100}}$ supported on integers all of whose prime factors are greater than $q_l$.  Then for $w\in [L_l, U_l]$ we have
		\begin{equation}
			\E_{+}\left[\left|Q_l(1/2,\chi)\right|^2 1(B_l\cap C_l\cap \left\{S_{n_l}\in (w,w+1]\right\}\right]\ll \E_{\oplus}\left[\left|Q_l(1/2,\chi)\right|^2\right]\frac{e^{-\frac{w^2}{n_l}}}{\sqrt{n_l}}.
		\end{equation}
	\end{lem}
	\begin{proof} 
		The proof is by approximating the indicator functions in terms of a short Dirichlet polynomial, as in \cite{AB} and \cite{ABR20}.
		We give the details for completeness.
		We define $\mathfrak{I}$ to be the set of $l$-tuples $\mathbf{u}=(u_1, \ ...,\ u_l)$ satisfying:
		\begin{equation}\label{I}        
			\mathlarger{\sum}^j_{i=1}u_i\in [L_j-1,U_j+1], \forall j< l, \qquad\mathlarger{\sum}^l_{i=1}u_i\in [w-1, w+2].
		\end{equation}
		Since $\kappa<1,$ we see for $\mathbf{u}$ to satisfy these bounds  on the partial sums $S_{n_j}$ for $1\le j\le l,$ we require for all $j>1$, \begin{equation}
			|u_j|\le 2(2\Delta_j+1).
		\end{equation}
		Moreover, we must have the inclusion
		\begin{equation}
			B_l\cap C_l\cap \{S_{n_l}\in (w,w+1]\}\subset \cup_{\mathbf{u}\in \mathfrak{I}}\{Y_j \in [u_j, u_j+\Delta_{j}^{-1}]\}.
		\end{equation}
		Hence, we obtain:
		\begin{equation}\label{expect}
			1(B_l\cap C_l\cap \{S_{n_l}\in (w,w+1]\}) \le \mathlarger{\sum}_{\mathbf{u}\in \mathfrak{I}}\prod_{1\le j\le l}1\{Y_j \in [u_j, u_j+\Delta_{j}^{-1}]\}.
		\end{equation}
		We then proceed to bound each indicator function with a short Dirichlet polynomial.
		
		Let $A\ge 10$, $j\le l$ and  $	\Delta_j=(n_j-n_{j-1})$.
		We take $D_{\Delta_{j},A}(x)$ to be the polynomial of Lemma 2.8 in \cite{AB}. Specifically, $D_{\Delta_{j},A}(x)$ is a polynomial of degree at most  $\Delta_{j}^{10A}$ such that
		\begin{equation}\label{indic}
			1(x \in [u_j+\Delta_j^{-1}])\le \left|D_{\Delta_{j},A}(x-u_j)\right|^2(1+ce^{-\Delta_{j}^{A-1}}),	
		\end{equation}
		where $c$ is an absolute constant. 
		If we let $X=Y_j$ then, by construction of the support of the coefficients of $Y_j,$ we see that  
		$D_{\Delta_j,A}(Y_j-u_j)$ is a Dirichlet polynomial on integers with all  prime factors in the interval $(q_{j-1},q_j]$ of length at most $\exp(2e^{n_j}\Delta_j^{10A})$.
		
		
		Hence we obtain
		\begin{equation}\label{approxindic}
			1(Y_j\in [u_j, u_j+\Delta_j^{-1}])\le \left|D_{\Delta_j,A}(Y_j-u_j)\right|^2(1+ce^{-\Delta_j^{A-1}}),
		\end{equation}
		where $Y_j=S_{n_j}-S_{n_{j-1}}$ for $1\le j\le l$.
		Putting this back in the expectation, we get 
		\begin{equation}\label{indictwist}
			\begin{aligned}
				&\E_{+}\left[\left|Q_l(1/2,\chi)\right|^2 1(B_l\cap C_l\cap \{S_{n_l}\in (w,w+1]\}\right]\\
				&\ll \mathlarger{\sum}\limits_{u\in \mathfrak{I}}			\E_{+}\left[\left|Q_l(1/2,\chi)\right|^2\prod_j(1+ce^{-\Delta_j^{A-1}})\left|D_{\Delta_j,A}(Y_j-u_j)\right|^2\right].
			\end{aligned}
		\end{equation}
		Following the decomposition of $|Q_l|^2$ into the contribution  from its real and imaginary parts as in the proof of Theorem \ref{twistthm}, we may apply Theorem \ref{realtwistthm} and Lemma \ref{Realsplit} to this combination of products of expectations of Dirichlet polynomials, $Q_j$, and real parts of Dirichlet polynomials, $D_{\Delta_j,A}(Y_j-u_j)$.
		Using Lemmas  \ref{Realsplit} and \ref{realprim} to remove the restriction on the characters being primitive  in the expectation in \eqref{indictwist}, and bounding the error terms $1+ce^{-\Delta_{j}^{A}-1}$ trivially, we see that the expectation is
		\begin{equation}\label{2.7link}
			\ll 	\E_{\oplus}\left[\left|Q_l(1/2,\chi)\right|^2\prod_j\left|D_{\Delta_j,A}(Y_j-u_j)\right|^2\right].
		\end{equation}
		
		Using Lemma \ref{Realsplit}, the expectation above splits, so the right-hand side is simply:
		\begin{equation}
			\E_{\oplus}\left[\left|Q_l(1/2,\chi)\right|^2\right]\prod_j\E_{\oplus}\left[\left|D_{\Delta_j,A}(Y_j-u_j)\right|^2\right].
		\end{equation}
		The orthogonality relation for the real part of even Dirichlet characters, Lemma \ref{Realexp}, implies that for each $1\le j\le l$,
		\begin{equation}
			\E_{\oplus}\left[\left|D_{\Delta_j,A}(Y_j-u_j)\right|^2\right]=\E\left[\left|D_{\Delta_{j},A}(\mathscr{Y}_j-u_j)\right|^2\right],
		\end{equation}
		where $(\mathscr{Y}_j,j\le l)$ are independent random variables of the form
		\begin{equation}
			\mathscr{Y}_j=\mathlarger{\sum}\limits_{q_{j-1}\le p \le q_j} \frac{\re X(p)}{p^{\frac{1}{2}}} .
		\end{equation}
		These are the exact analogue of the random variables defined in Equation 54 in \cite{AB}. Their distribution doesn't depend on being defined for the $t$-aspect or the $q$-aspect.
		These variables can be compared to Gaussian explicitly as in \cite{AB} starting at Equation 54. Proceeding as such, we obtain the result.
	\end{proof}  
	We also require a lemma which is similar to Lemma \ref{2.6}, but includes the mollifier $LM_1....M_l$ in the moment. This is analogous to Lemma 2.7 in \cite{AB}.
	\begin{lem}\label{2.7} Let $1\leq l\leq \mathcal L$. Let $Q_l$ be a Dirichlet polynomial as defined in Lemma \ref{2.6}, and $w\in [L_l,U_l]$.
		Then 
		\begin{equation}\label{twistindic}
			\E_{+}\left[\left|LM_1,...M_{l}Q_l\right|^21(B_l\cap C_l, S_{n_l}\in (w,w+1])\right]\ll   \E_{\oplus}\left[\left|Q\right|^2\right]\times\frac{\log q}{\log q_l} \frac{e^{-\frac{w^2}{n_l}}}{\sqrt{n_l}}.
		\end{equation}
	\end{lem}
	\begin{proof}
		The proof follows from the proof of Lemma \ref{2.6}. We see using equations \eqref{indic} and \eqref{approxindic}, the left-hand side of \eqref{twistindic} becomes 
		\begin{equation}\label{twistsplit}
			\ll (1+ce^{-\Delta_j^{A-1}})\mathlarger{\sum}\limits_{u\in \mathfrak{I}} \E_{\oplus} \left[\left|LM_1...M_l\right|^2\left|Q_l\right|^2\prod_j \left|D_{\Delta_j,A}(Y_j-u_j)\right|^2 \right].
		\end{equation}
		Setting \begin{equation}
			Q=Q_l\prod_j D_{\Delta_{j},A}(Y_j-u_j),
		\end{equation}
		we see that $Q$ is a well-factorable Dirichlet polynomial.
		Applying Theorem \ref{twistthm} shows that the expression in \eqref{twistsplit} is 
		\begin{equation}
			\ll (1+ce^{-\Delta_j^{A-1}})\frac{\log q}{\log q_l}\times  \mathlarger{\sum}\limits_{u\in \mathfrak{I}} \E_{\oplus} \left[\left|Q_l\right|^2\prod_j \left|D_{\Delta_j,A}(Y_j-u_j)\right|^2 \right].
		\end{equation} 
		Proceeding as from Equation \ref{2.7link}, we obtain Equation \ref{twistsplit}. 
	\end{proof}
	We begin by showing that for $l\ge 1,$
	\begin{equation}\P(A_{l+1}^c\cap G_l)\ll \frac{e^{-\frac{V^2}{\log\log q}}}{\sqrt{\log\log q}}(\log_{l+1}(q))^{-\delta}.\end{equation}
	For any $k>1$, we have, $\P(A_{l+1}^c\cap G_l)$ is bounded by: \begin{equation}\label{notA}\mathlarger{\sum}\limits_{u\in [L_l, U_l]}\E_{+}\left[\frac{|\tilde{S}_{n_{l+1}}-\tilde{S}_{n_l}|^{2k}}{(\mathscr{A}(n_{l+1}-n_l))^{2k}}1(B_l\cap C_l\cap \{S_{n_l}\in (u,u+1]\})\right].\end{equation}
	We let $k=\lceil\mathscr{A}^2(n_{l+1}-n_l).\rceil$
	The Dirichlet polynomial $\left(\tilde{S}_{n_{l+1}}-\tilde{S}_{n_l}\right)^k$ satisfies the conditions of Lemma \ref{2.6} for each $u\in [L_l,U_l],$ so we may bound the expression in Equation \eqref{notA} as:
	\begin{equation}
		\ll \frac{\E_{\oplus}\left[|\tilde{S}_{n_{l+1}}-\tilde{S}_{n_l}|^{2k}\right]}{(\mathscr{A}(n_{l+1}-n_l))^{2k}}\mathlarger{\sum}\limits_{u\in [L_l,U_l]}\frac{e^{-\frac{u^2}{n_l}}}{\sqrt{n_l}}.  \end{equation}
	We can then use Lemma \ref{Strongmoment} to bound the moment, giving the above expression as:
	\begin{equation}
		\ll \frac{k! (n_{l+1}-n_l+1)^k}{(\mathscr{A}(n_{l+1}-n_l))^{2k}}\mathlarger{\sum}\limits_{u\in [L_l,U_l]}\frac{e^{-\frac{u^2}{n_l}}}{\sqrt{n_l}}.  
	\end{equation}
	
	We apply Stirling's formula to bound the above expression as:
	\begin{equation}
		\ll \frac{\sqrt{k} \left(\frac{\mathscr{A}^2(n_{l+1}-n_l+1)^2}{e}\right)^k}{(\mathscr{A}(n_{l+1}-n_l))^{2k}}\mathlarger{\sum}\limits_{u\in [L_l,U_l]}\frac{e^{-\frac{u^2}{n_l}}}{\sqrt{n_l}} 
		\ll e^{-\mathscr{A}^2 (n_{l+1}-n_l)}\sqrt{n_{l+1}-n_l} \mathlarger{\sum}\limits_{r\ge -\mathbf{s}\log_{l+2}q}\frac{e^{\frac{-(\kappa n_l+r)^2}{n_l}}}{\sqrt{n_l}} .
	\end{equation}
	Since $L_l=\kappa n_l-\mathbf{s} \log_{l+2} q$, the above expression is
	\begin{equation}
		\begin{aligned}
			\label{bound3}
			\ll  \frac{\ e^{-\kappa^2 n_l }}{\sqrt{n_l}} (\log_{l+1}q)^{-\mathscr{A}^2\mathbf{s}+2\kappa\mathbf{s}} (\log_{l+2}q)^{\frac{1}{2}+\mathscr{A}^2s}
			\ll \frac{e^{-\kappa^2\log\log q}}{\sqrt{\log\log q}} (\log_{l+1} q)^{\kappa^2s-\mathscr{A}^2\mathbf{s}+2\kappa\mathbf{s}+1}.
		\end{aligned}
	\end{equation}
	The choice of parameters for Equation \ref{Afix} ensures the exponent is negative.
	
	Next, we bound $\P(B^c_{l+1}\cap G_l).$ Here again we have a sharper bound that the analogous result in \cite{AB}, which allows us to take $V=o(\log\log q).$
	We have: \begin{equation}
		\P(B^c_{l+1}\cap G_l)\ll \sum\limits_{u\ge U_{l+1}}\P (\left\{S_{n_{l+1}}\in [u,u+1]\right\}\cap G_l).\end{equation}
	Using Lemma \ref{2.6} with $Q_l=1$, we see this is:
	\begin{equation}
		\label{geoprog}
		\ll\sum\limits_{u\ge U_{l+1}} \frac{e^{-\frac{u^2}{n_{l+1}}}}{\sqrt{\log \log q_l}}\ll \frac{1}{\sqrt{\log \log q}} e^{-\frac{U_{l+1}^2}{n_{l+1}}}\sum\limits_{k\ge 0} e^{\frac{U_{l+1}^2-(U_{l+1}-k)^2}{n_{l+1}}}.\end{equation}
	by writing $k=\lfloor U_{l+1}-u\rfloor$.
	The sum over $k$ is $\ll\sum\limits_{k\ge 0} e^{-\frac{2k U_{l+1}}{n_{l+1}}}\ll1$.
	So we see that \eqref{geoprog} is
	\begin{equation}
		\begin{aligned}
			\ll \frac{e^{-\frac{U_{l+1}^2}{n_{l+1}}}}{\sqrt{\log \log q}} 
			\ll\frac{e^{-\kappa^2 \log\log q}}{\sqrt{\log\log q}} e^{\kappa^2 \mathbf{s}\log_{l+3} q-2\kappa \mathbf{s} \log_{l+3}q}
			\ll \frac{e^{-\kappa^2 \log\log q}}{\sqrt{\log\log q}} (\log_{l+2} q)^{\kappa (\kappa-2)\mathbf{s}}.
		\end{aligned}
	\end{equation}
	The choice of parameters ensures the exponent is negative.
	
	We now turn to bound $\P(H\cap C_{l+1}^c \cap A_{l+1}\cap D_{l+1}\cap G_l)$. By using Lemma \ref{23ABR} and partitioning on the values $S_{n_l}=u$ according to the range allowed by the event $G_l$, as well as the values $v=S_{n_{l+1}}-S_{n_l}$  allowed by the event $A_{l+1}\cap C_{l+1}^c,$ we see that the probability may be bounded by:
	\begin{equation}\label{impcons} 
		\mathlarger{\sum}\limits_{\substack{u\in [L_l,U_l]\\ u+v\le L_{l+1}\\
				|v|\le \mathscr{A}(n_{l+1}-n_l)}} \P(\{S_{n_l}\in (u,u+1], S_{n_{l+1}}-S_{n_l}\in (v,v+1], \left|LM_1...M_{l+1}\right|>\frac{1}{100} e^{V-(u+v)}\}\cap B_l\cap C_l).\end{equation}
	We improve on the analogous bound in \cite{AB} by using \eqref{hardS}, which is a better upper bound on the moment of the truncated sums $S_{n_{l+1}}-S_{n_l}$ than was achieved in Equation (79) in \cite{AB}. This allows us to attain an implicit constant in Theorem \ref{LargeDev} which is bounded as $\alpha\to 0,$ and gives us the bounds to prove Theorem \ref{largerange} for values of $V$ where we may have $V=o(\log\log q).$
	
	The probability in \eqref{impcons} may be bounded  by:
	\begin{align} 
		\mathlarger{\sum}_{\substack{u\in [L_l,U_l]\\ u+v\le L_{l+1}\\
				|v|\le \mathscr{A}(n_{l+1}-n_l)}} e^{-2(u+v-V)} \E_+[ |LM_1...M_{l+1}|^21(\{S_{n_l}\in (u,u+1], S_{n_{l+1}}-S_{n_l}\in (v,v+1]\}\cap B_l\cap C_l)].
	\end{align}	
	By Markov's inequality, if we set $r=\left\lceil \frac{v^2}{(n_{l+1}-n_l)}\right\rceil$, then this is:
	\begin{equation}
		\ll\sum\limits_{\substack{u\in [L_l,U_l]\\
				u+v\le L_{l+1}\\	|v|\le \mathscr{A}(n_{l+1}-n_l)}}e^{-2V+2(u+v)}\E\left[ \frac{\left| S_{n_{l+1}}-S_{n_l}\right|^{2r}}{v^{2r}} \left|LM_1.....M_{l+1}\right|^21(\left\{S_{n_l}\in (u,u+1]\right\}\cap B_l\cap C_l)\right].
	\end{equation}
	Using Lemma \ref{2.7} with $Q_l=(S_{n_{l+1}}-S_{n_l})^r$ and \eqref{hardS} yields that this is
	\begin{equation}
		\ll\sum\limits_{\substack{
				u+v\le L_{l+1}}}e^{-2V+2(u+v)}  e^{2(\log\log q-n_{l+1})}\frac{(2r)!}{(2v)^{2r} r!}(n_{l+1}-n_l)^r\times \frac{e^{-\frac{u^2}{n_l}}}{\sqrt{n_l}}.\end{equation}
	We can use Stirling's formula to bound this as 
	\begin{equation}
		\ll\frac{1}{\sqrt{\log \log q}}\sum\limits_{\substack{
				u+v\le L_{l+1}}}e^{-2V+2(u+v)}  e^{2(\log\log q-n_{l+1}})(\frac{r(n_{l+1}-n_l)}{ev^2})^r e^{-\frac{u^2}{n_l}}.\end{equation}
	Substituting the value for $r$, this is:
	\begin{equation}
		\ll\frac{1}{\sqrt{\log \log q}}\sum\limits_{\substack{
				u+v\le L_{l+1}}}e^{-2V+2(u+v)}  e^{2(\log\log q-n_{l+1})}e^{-\frac{v^2}{n_{l+1}-n_l}}e^{-\frac{u^2}{n_l}}.\end{equation}
	
	Making the change of variables $a=u-\kappa n_l, b=v-\kappa(n_{l+1}-n_l)$ shows this is:
	\begin{equation}\label{Gaussum}
		\ll\sum\limits_{a+b \le L_{l+1}-\kappa n_{l+1}}\frac{e^{-\kappa^2 n_{l+1}}}{\sqrt{\log \log q}} e^{2 (\log \log q-n_{l+1})}e^{(2-2\kappa)(a+b)}e^{-\frac{b^2}{n_{l+1}-n_l}}.	
	\end{equation}
	Performing the sum over $a$ gives that \eqref{Gaussum} is 
	\begin{equation}
		\ll \sum\limits_b \frac{e^{-\kappa^2 n_{l+1}}}{\sqrt{\log \log q}} e^{2 (\log \log q-n_{l+1})} e^{-\frac{b^2}{n_{l+1}-n_l}-(2-2\kappa)\mathbf{s} \log_{l+2}q}.
	\end{equation}
	Using a Gaussian bound, the above expression is:
	\begin{equation}\label{bound4}
		\ll	\frac{e^{-\kappa^2 n_{l+1}}}{\sqrt{\log \log q}} e^{2 (\log \log q-n_{l+1})} \sqrt{n_{l+1}-n_l}e^{-(2-2\kappa)\mathbf{s} \log_{l+2}q}\\ 
		\ll\frac{e^{-\kappa^2 \log\log q}}{\sqrt{\log \log q}} (\log_{l+1} q)^{\frac{1}{2}+\kappa^2-(2-2\kappa)\mathbf{s}}.\end{equation}
	By construction of $\mathbf{s}$, the exponent of $\log_{l+1} q$ is negative. We remark that $\log_{l+1}q\gg \log_{l+2} q,$ so that this bound is even sharper than the total bound for Proposition \ref{prop G}(ii).
	
	It remains to bound $\P(D_{l+1}^c\cap A_{l+1}\cap G_l).$	The event $D_{l+1}^c\cap A_{l+1}\cap G_l$ is contained in \begin{equation}
		\{\left|LM_1...M_l\right|>e^{\mathscr{A}(\log\log q -n_l)}\}\cap G_l.
	\end{equation}
	This is because $\{\left|LM_1...M_l\right|\leq e^{\mathscr{A}(\log\log q -n_l)}\}\cap D_l\}$ is in $A_{l+1}\cap D_{l+1}$ by using Lemma \ref{23ABR}.
	By Markov's inequality, we have
	\begin{equation}
		\P(\{\left|LM_1...M_l\right|>e^{\mathscr{A}(\log\log q -n_l)}\}\cap G_l)\ll  e^{-2\mathscr{A}(\log\log q-n_l)}\E_{\oplus}\left[\left|LM_1...M_l\right|^2 1(G_l)\right].
	\end{equation}
	Using Lemma \ref{2.7}, we can bound the right-hand side as  
	\begin{equation}
		\begin{aligned}
			&\ll   e^{-2(\mathscr{A}-1)(\log\log q-n_l)}\times \frac{e^{-\frac{L_l^2}{n_l}}}{\sqrt{n_l}}
			&\ll \frac{e^{-\kappa^2 \log\log q}}{\sqrt{\log\log q}}(\log_{l+1}q)^{2\mathbf{s}\left(\kappa+1 -\mathscr{A}\right)}.
		\end{aligned}
	\end{equation}
	The choice of $\mathscr{A}$ ensures the exponent is negative.
	
	\subsection{Proof of Proposition \ref{prop G}(iii)}
	We have the obvious bound $\PP(H\cap G_{\mathcal L})\leq \PP(G_{\mathcal L})$. 
	It then suffices to apply Lemma \ref{2.7}. Note that the length of the interval $[L_{\mathcal L}, U_{\mathcal L}]$ is of order one. 
	Therefore,
	\begin{equation}
		\PP(G_{\mathcal L})\ll \sup_{v\in [L_{\mathcal L}, U_{\mathcal L}]} 
		\frac{1}{\sqrt{n_{\mathcal L}}}e^{-v^2/{n_{\mathcal L}}}\ll \frac{e^{-\kappa^2\log\log q}}{\sqrt{\log\log q}}.
	\end{equation}
	
	%
	%
	%
	
	\section{Proof of Corollary \ref{momentcor}}\label{CorSec}
	Corollary \ref{momentcor} follows from Theorem \ref{LargeDev} in a similar way to the proof of Corollary 1.2 in \cite{AB}.
	\begin{proof} Consider the CDF of the random variable $\log L(1/2,\chi)$ with $\chi$ drawn uniformly from the even primitive characters mod $q$, i.e., $F(V)=\P(\log\left|L(1/2,\chi)\right|\le V)$.
		Then if we set \begin{equation}
			\mathscr{M}_\beta=	\E_{+}\left[ \left|L(1/2,\chi)\right|^{2\beta}\right] ,
		\end{equation}
		for $0<\beta < 1$,	the moments can be written as: 
		\begin{equation}
			\label{moment}
			\mathscr{M}_\beta=\int^{+\infty}_{-\infty} e^{2\beta V}d F(V)=\left[e^{2\beta V}S(V)\right]^{+\infty}_{-\infty}+\int^\infty_{-\infty}2\beta e^{2\beta V}S(V)dV.
		\end{equation}
		Since $S(V)\le 1$, the boundary term at $-\infty$ is zero.
		Since there are only finitely many characters, the boundary term at $+\infty$ is also zero.	The contribution of negative values of $V$ in \eqref{moment} is negligible, since \begin{equation}\int^0_{-\infty}2\beta e^{2\beta V} S(V)dV \le \int^0_{-\infty}2\beta e^{2\beta V}dV=1.\end{equation}
		
		In order to estimate $\int^\infty_{0}2 \beta e^{2\beta V}S(V)$, consider $\beta_-$ and $\beta_+$ with \begin{equation}0<\beta_-<\beta<\beta_+<1.\end{equation}
		We want to show the dominant contribution of to the integral in \eqref{moment} comes from the interval $\left[\beta_-\log\log q, \beta_+\log\log q\right]$.
		We take 
		$\beta_-=\frac{\beta}{4}, \quad \beta_+=\frac{3+\beta}{4}.$
		Using Theorem \ref{LargeDev} to bound $S(V),$ we have
		\begin{equation}
			\begin{aligned}
				\int^{\beta_+\log\log q}_{\beta_-\log\log q}\beta e^{2\beta V}S(V)dV&
				\ll \int^{\beta_+\log\log q}_{\beta_-\log\log q} e^{2\beta V}\frac{e^{-\frac{V^2}{\log\log q}}}{\sqrt{\log\log q}}dV\\
				&=\frac{e^{\beta^2\log\log q}}{{\sqrt{\log\log q}}}\int^{\beta_+\log\log q}_{\beta_-\log\log q}e^{-(\beta\log\log q-V)^2/\log\log q}dV\\&\ll e^{\beta^2\log\log q}.
			\end{aligned}
		\end{equation}
		The contribution of the interval $[0, \beta_-\log\log q]$ to the integral in \eqref{moment} is of a lower order since it is bounded by:
		\begin{equation}
			\int^{\beta_-\log\log q}_02\beta e^{2\beta V}dV \le e^{\frac{\beta^2\log\log q}{2}}.
		\end{equation}
		
		For the interval $\left[\beta_+\log\log q, \infty\right),$ we require a bound on $S(V) $ in this range. But we have \begin{equation}\E_{+}[\left|L(1/2,\chi)\right|^2]\ll \log q,\end{equation}
		and so Markov's inequality yields $	S(V)\ll e^{-2V}\log q$.
		Hence the contribution in the range $[\beta_+\log\log q, \infty)$ to the integral in \eqref{moment} is 
		\begin{equation}\ll{\log q} \int^\infty_{\beta_+\log\log q} 2\beta e^{2(\beta-1)V} dV
			\label{maincon}= \frac{\beta}{1-\beta}e^{(2(\beta -1)\beta_++1)\log\log q}.\end{equation}
		This is $o(e^{\beta^2\log\log q})$ by the choice of  $\beta_+.$
	\end{proof}
	Since we chose the parameters in Section \ref{choiceparameters}  such that the inequalities held uniformly for $0<\kappa<1$, we see that the implicit constant in Equation \ref{momentcor} is bounded as $\beta\to 0$. As previously stated, this contrasts with the proof in \cite{AB}, where in Equation \ref{blowup}, their bound had $C_\beta\to \infty$ as $\beta\to 0$.

	\bibliographystyle{alpha}
	\bibliography{ref.bib}

\begin{thebibliography}{BPRZ20}

\bibitem[AB23]{AB}
L.-P. Arguin and E.~Bailey.
\newblock Large deviation estimates of {S}elberg's {C}entral {L}imit {T}heorem
  and applications.
\newblock {\em Int. Math. Res. Not. IMRN}, (23):20574--20612, 2023.

\bibitem[ABR20]{ABR20}
L.-P. {Arguin}, P.~{Bourgade}, and M.~{Radziwi{\l}{\l}}.
\newblock {The Fyodorov-Hiary-Keating Conjecture. I}.
\newblock {\em arXiv e-prints}, page arXiv:2007.00988, July 2020.

\bibitem[BPRZ20]{Bui}
H.~M. Bui, K.~Pratt, N.~Robles, and A.~Zaharescu.
\newblock Breaking the {$\frac12$}-barrier for the twisted second moment of
  {D}irichlet {$L$}-functions.
\newblock {\em Adv. Math.}, 370:107175, 40, 2020.

\bibitem[Gao21]{gao}
P.~Gao.
\newblock {Bounds for moments of Dirichlet $L$-functions to a fixed modulus}.
\newblock {\em arXiv e-prints}, page arXiv:2103.00149, February 2021.

\bibitem[HRS19]{HRS19}
W.~Heap, M~Radziwi\l{\l}, and K.~Soundararajan.
\newblock {Sharp upper bounds for fractional moments of the {R}iemann Zeta
  Function}.
\newblock {\em Q. J. Math.}, 70(4):1387--1396, 2019.

\bibitem[HW20]{PW}
P.-H. Hsu and P.-J. Wong.
\newblock On {S}elberg's central limit theorem for {D}irichlet {$L$}-functions.
\newblock {\em J. Th\'{e}or. Nombres Bordeaux}, 32(3):685--710, 2020.

\bibitem[IS99]{Sarnak}
H.~Iwaniec and P.~Sarnak.
\newblock {Dirichlet $L$-functions at the central point}.
\newblock {\em Number theory in progress}, 2:941--952, 1999.

\bibitem[KS00]{KS}
J.~P. Keating and N.~C. Snaith.
\newblock Random matrix theory and {$\zeta(1/2+it)$}.
\newblock {\em Comm. Math. Phys.}, 214(1):57--89, 2000.

\bibitem[Pra19]{Pratt}
K.~Pratt.
\newblock Average nonvanishing of {D}irichlet {$L$}-functions at the central
  point.
\newblock {\em Algebra Number Theory}, 13(1):227--249, 2019.

\bibitem[RS05]{RS}
Z.~Rudnick and K.~Soundararajan.
\newblock Lower bounds for moments of {$L$}-functions.
\newblock {\em Proc. Natl. Acad. Sci. USA}, 102(19):6837--6838, 2005.

\bibitem[Sel46]{Sel}
A.~Selberg.
\newblock Contributions to the theory of the {R}iemann zeta-function.
\newblock {\em Arch. Math. Naturvid.}, 48(5):89--155, 1946.

\bibitem[Sou09]{Sound}
K.~Soundararajan.
\newblock Moments of the {R}iemann zeta function.
\newblock {\em Ann. of Math. (2)}, 170(2):981--993, 2009.

\end{thebibliography}

\end{document}